\documentclass[final]{amsart} 

\usepackage[utf8]{inputenc}
\usepackage[T1]{fontenc}
\usepackage{lmodern}

\usepackage{graphicx}
\usepackage{subcaption}
\usepackage{paralist}

\usepackage{tikz}
\usetikzlibrary{external, arrows.meta, positioning, shadings, math}
\usepackage{pgfplots}
\pgfplotsset{compat=1.17}
\pgfplotscreateplotcyclelist{my color list}{
    {blue!40!black,mark=none},
    {red!40!black,mark=none},
    {green!40!black,mark=none},
    {blue!40!black,dashed},
    {red!40!black,dashed},
    {green!40!black,dashed}
}

\usepackage{mathtools}
\usepackage{amsmath}
\usepackage{amssymb}
\usepackage{algorithm}
\usepackage{algpseudocode}
\usepackage{algorithmicx}

\usepackage{comment}

\usepackage{todonotes}
\usepackage{hyperref}
\hypersetup{
    citecolor=green,
    colorlinks=true,
    linkcolor=blue,
    filecolor=magenta,      
    urlcolor=red,
}

\usepackage[maxnames=10, bibencoding=utf8]{biblatex}
\renewcommand{\div}{\mathrm{div}}
\newcommand{\id}{\mathbf{id}}

\title[$p$-Harmonic Descent Approach in Fluid Dynamic Shape Optimization]{A novel $p$-Harmonic Descent Approach applied to Fluid Dynamic Shape Optimization}

\author[P. M. M{\"u}ller, N. K{\"u}hl, M. Siebenborn, K. Deckelnick, M. Hinze, T. Rung]{Peter Marvin M{\"u}ller, Niklas K{\"u}hl, Martin Siebenborn,\\Klaus Deckelnick, Michael Hinze, Thomas Rung}

\date{March 26, 2021}

\begin{document}
	
\begin{abstract}
We introduce a novel method for the implementation of shape optimization in fluid dynamics applications, where we propose to use the shape derivative to determine deformation fields with the help of the $p-$ Laplacian for $p > 2$. This approach is closely related to the computation of steepest descent directions of the shape functional in the $W^{1,\infty}-$ topology and refers to the recent publication \cite{DHH21}, where this idea is proposed. Our approach is demonstrated for shape optimization related to drag-minimal free floating bodies. The method is validated against existing approaches with respect to convergence of the optimization algorithm, the obtained shape, and regarding the quality of the computational grid after large deformations. Our numerical results strongly indicate that shape optimization related to the $W^{1,\infty}$-topology -- though numerically more demanding -- seems to be superior over the classical approaches invoking Hilbert space methods, concerning the convergence, the obtained shapes and the mesh quality after large deformations, in particular when the optimal shape features sharp corners.
\end{abstract}

\maketitle

\section{Introduction}
\label{sc:introduction}
Adjoint-based local optimization has been matured towards an efficient industrially applied strategy, e.g. \cite{othmer2014adjoint, papoutsis2016continuous}. When attention is given to shape optimization, the aim is to find optimal shapes regarding a physical quantity $J$, e.g. the drag force experienced by an obstacle. Mathematically speaking, a shape functional $J$ is minimized subject to partial differential equation (PDE) constraints. The latter typically govern the physics, e.g. the conservation of mass and momentum. A crucial part of the shape optimization procedure is the choice of the descent direction that provides an  update rule for the design variable. In our case the design variable refers to the underlying geometry $\Omega$, the boundary $\partial \Omega$ or parts of the boundary $\Gamma$. The descent direction is usually employed to converge a sequence of shape updates via a deformation field $\mathbf{u}$. First attempts used the directional derivative $J'$, also referred to as shape derivative  in order to formulate a shape update rule \cite{sokolowski1992, kawohl1998optimal, zolesio2011}. Because shapes are not elements of a vector space, e.g. there is no meaningful definition of the summation of two shapes, the shape derivative is defined by introducing deformations that make shapes variable and thus allows the definition of directional derivatives. With the help of shape calculus, cf. \cite{sokolowski1992, zolesio2011} for a mathematical perspective or \cite{schmidt2013three, kuhl2019decoupling} for an engineering perspective, a shape derivative can be computed that relates to a scalar field $\gamma$ defined on the boundary $\Gamma$. We want to refer to $\gamma$ as the local shape sensitivity. Mind that the computation of the shape derivative with shape calculus is a rather involved task and heavily depends on the objective functional $J$, the PDE constraints that apply as well as the domain $\Omega$ and the boundary $\partial \Omega$, respectively. In \cite{pironneau1973optimum} the shape is updated by using the local shape sensitivity $\gamma$ to perform the deformation in normal direction $\mathbf{n}$ of the boundary $\Gamma$. Using the descent direction $\mathbf{u} = -\gamma \mathbf{n}$ yields a contribution $\mathbf{u} \cdot \mathbf{n} = - \gamma^2$ to reducing the objective functional, and is a popular approach within shape optimization, cf. \cite{vassberg2006aerodynamic, vassberg2006aerodynamic_II}. The attempt is limited since it often yields shapes with rough/noisy boundaries \cite{stuck2011adjoint, kroger2015cad} and distorted near-wall meshes which in turn hamper the preservation of numerical accuracy during the optimization procedure \cite{stavropoulou2014plane, bletzinger2014consistent}. In fluid dynamics and neighboring applications, the computational grid is frequently morphed and not renewed after each optimization step. Thereupon, attempts follow in order to gain higher regularity of the deformation. One approach is based on the definition of a shape gradient $\mathrm{grad}J$ by an inner product and the shape derivative $J'$. Here the shape gradient is identified by the Riesz representation of the directional derivative of the shape functional. Even though this leads to smoother deformations, the approach is algorithmic challenging due to solving a PDE on a hyperplane, and also mathematically questionable in the general case, see also \cite{allaire2020} and cf. section~\ref{sc:framework}. In this regard, different gradients are associated with different transformations applied to the shape derivative, and several techniques have been proposed to increase the regularity of the shape updates:
\begin{itemize}
\item [a)] CAD-related shape definitions connect the node-based shape derivatives to the CAD parameterization using the chain rule of differentiation, cf. \cite{lohner2003adjoint, robinson2012optimizing}. The procedure couples the various local derivatives and thereby ensures smooth shapes. However, the rigid finite dimensional initial CAD parameterization limits the attainable shapes  and different CAD models may result in different optimal shapes.
\item [b)] A coupling of mesh node updates using either local shape functions, e.g. FE-type  functions \cite{soto2002stabilized, soto2004adjoint}, or global shape functions, e.g. Hicks-Henne approaches \cite{hicks1978wing}.
\item [c)] A more rigor approach of Jameson and Vassberg \cite{jameson2000studies, vassberg2006aerodynamic, vassberg2006aerodynamic_II} applies an implicit, continuous smoothing operator to either the shape derivative $J'$ or the deformation field $\mathbf{u}$, based on an extended definition of the inner product, frequently labeled 'Sobolev-gradient'. Applying the smoothing operation on the surface leads to the Laplace-Beltra\-mi operator and a related surface metric. For computational reasons, the practice is often performed in an explicit manner, cf. \cite{bletzinger2014consistent}. The explicitly filtered local shape sensitivity, e.g. by using consistent kernel functions \cite{kroger2015cad}, marks a first-order approximation to the implicit Sobolev-gradient \cite{stuck2011adjoint}.
\end{itemize}
All strategies (a-c) essentially couple node updates and thereby obtain smooth design updates. Having updated the discrete design surface, the subsequent numerical investigation of the updated design also requires an update of the computational mesh, i.e. the shape gradient of the design surface needs to be extended into the domain.  The habitat of the shape gradient depends on the surface metric and can be surface as well as volume based. Prominent examples refer to above mentioned Laplace-Beltrami (LB) or the Steklov-Poincar{\'e} (SP) metric, c.f. \cite{schulz2016computational}. The first approach (LB) exclusively operates in the tangent space of the design surface, the latter (SP) leads to a domain formulation, where results are subsequently projected on the controlled shape. The SP strategy gives the shape update of the design surface and mesh using the shape sensitivity \cite{schulz2016computational, haubner2020continuous, allaire2020}. The volume based SP approach is particularly attractive for optimization procedures which prefer mesh morphing over re-meshing strategies, as it is customary for engineering simulations. While re-meshing can be automated, the lack of fair restart capabilities  becomes prohibitively expensive in practical applications. Moreover, the use of standardized, HPC capable solution routines supplied by the flow solver (assembling, solving, etc.) represent another significant benefit of the SP approach.

A different avenue is taken by the phase field method for fluid mechanic shape optimization, where the shape of the sought domain is approximated by the zero level set of a phase field function. This turns the shape optimization problem into a PDE constrained optimization problem where the phase field enters as control in the coefficients of the PDE. This allows to apply the complete algorithmic machinery of PDE constrained optimization methods to this formulation of the fluid dynamic shape optimization problems, and this approach also naturally allows topology changes of the shape. However, numerical methods also for this approach encounter problems in situations where the sought shape needs to develop kinks and/or corners. This approach is proposed in \cite{BP03}, and in a couple of papers investigated for hydrodynamic shape optimization problems \cite{GHHK15,GHHKL16,GHKL18}.

Although smooth shapes may be desirable for different reasons, they are not necessarily optimal. If the optimum involves a kink or a corner, the above mentioned strategies to obtain smooth shape updates display difficulties to capture such optima from curved initial configurations, if the respective region is  (initially) not resolved by very fine grids. Section \ref{sc:numerics} of the present paper discusses a classical example of a pointed optimal shape. The same is true for counterpart situations, i.e. to transform an initially kinked shape into a curved optimum. Though this might be possible, the convergence is often fairly slow.

The present study aims to convey the merits of an alternative strategy 
to compute the shape deformation from shape derivatives in the context of CAD-free -- aka node-based -- shape optimization.
The $p$-Laplace operator is used in a volume-based formulation along the route of the SP metric. The approach is industrially feasible and supports unstructured meshes. Applications refer to $2D$ and $3D$ fluid dynamic shape optimization ranging from laminar to turbulent external flows around free floating objects with fixed displacement.

The remainder of the paper is structured as follows: Section~\ref{sc:framework} outlines the mathematical framework and the rationale that leads to the $p$-Laplace problem to approximate the steepest descent direction within the $W^{1,\infty}$ - topology together with a discussion of our fluid dynamic shape optimization problem. Section~\ref{sc:algorithm} presents the solution algorithm. Section~\ref{sc:numerics} applies the approach to three test cases and the manuscript closes with conclusions in Sec. \ref{sc:conclusion}.

\section{Mathematical Framework}
\label{sc:framework}
In this section we want to outline the basic idea behind the $p$-harmonic approach and briefly recall the concept of shape optimization. 

For this purpose let $J: \mathcal A \rightarrow \mathbb{R}$ denote a shape functional, where $\mathcal A$ denotes the set of admissible domains which has to be specified in the respective application. In our setting the set $\mathcal A$ is specified through Fig. \ref{fig:domain}, and the shape functional with \eqref{Forcefunctional} is given in \eqref{eq:optProblemObjectiveFunction}. For the algorithmic minimization of $J$ for a given domain $\Omega \in \mathcal A$ we intend to specify descent vector fields $\mathbf{u}^\ast: \mathbb{R}^d \rightarrow \mathbb{R}^d$ such that $J'(\Omega;\mathbf{u}^*) <0$ holds, where $J'(\Omega;\mathbf{u}^*)$ denotes the shape derivative of $J$ at $\Omega$ in direction $\mathbf{u}^*$. The perturbed domain $\tilde\Omega$ then has the form
\[
\tilde\Omega = \mathbf{T}_t(\Omega) := (\id + t \mathbf{u}^*)(\Omega),
\]
where $t>0$ is a step size specified in the respective minimization algorithm. However, descent in this context requires to specify an appropriate topology. Moreover, it frequently is required that Lipschitz domains $\Omega$ are mapped to Lipschitz domains $\tilde\Omega$. Common practical approaches use Hilbert Space methods, i.e. seek descent vector fields $\mathbf{u}^*$ determined with the help of the shape derivative by
\[
a(\mathbf{u}^*, \mathbf{w}) = J'(\Omega; \mathbf{w}) \text{ for all } \mathbf{w} \in H,
\]
where $(H,a(\cdot,\cdot))$ denotes an appropriate Hilbert space, see e.g. \cite[Section 5-9]{allaire2020} for an extensive discussion of approaches related to the Hilbert space setting. This requires to compute the Riesz representative of the functional $J'(\Omega;\cdot)$. A typical choice of $H$ is $H^m(\Omega,\mathbb{R}^d)$, where however $m \in \mathbb{N}$ has to be chosen large enough to obtain a Lipschitz transformation. A way around this would be to directly choose $\mathbf{u}^* \in W^{1,\infty}(\Omega,\mathbb{R}^d)$ as a direction of steepest descent for $J$ at $\Omega$, where $W^{1,\infty}(\Omega,\mathbb{R}^d)$ denotes the set of Lipschitz transformations from $\Omega$ to $\mathbb{R}^d$. This leads to the minimization problem
\begin{equation}
    \min_{\mathbf{u} \in W^{1,\infty}(\mathbb{R}^d,\mathbb{R}^d), \, \| \mathbf{u} \|_{W^{1,\infty}} \le 1} \quad J'(\Omega; \mathbf{u}),
    \label{eq:infinityLaplacianProblem}
\end{equation}
which is challenging both from the mathematical and the numerical perspective, since it represents a variational problem in the non-reflexive Banach space $W^{1,\infty}(\Omega,\mathbb{R}^d)$. Variational problems of this kind are studied in e.g. \cite{ishii2005limits}, where it  is proposed to approach solutions $\mathbf{u}^*$ of problem \eqref{eq:infinityLaplacianProblem} by a sequence of solutions $\mathbf{u}_p^* \in W^{1,p}(\Omega,\mathbb{R}^d)$ of the variational problem
\begin{equation}
    \min_{u_p\in W^{1,p}(\Omega, \mathbb{R}^d)} \; \frac{1}{p} \int_\Omega (\nabla \mathbf{u}_p \colon \nabla \mathbf{u}_p)^{\frac{p}{2}} \, \mathrm{d} x
    + J'(\Omega; \mathbf{u}_p) 
    \label{eq:pLaplacianVariationalForm}
\end{equation}
for $p>2$, compare \cite[Proposition 5.2,5.3]{ishii2005limits} for the mathematical analysis of the limit process $p^* \rightarrow \infty$. Note that for $p=2$ we recover a Hilbert space setting as described above, and refer to problem \eqref{eq:pLaplacianVariationalForm} as $p-$Laplace relaxation of problem \eqref{eq:infinityLaplacianProblem}. In the next section we adapt problem \eqref{eq:pLaplacianVariationalForm} to our fluid dynamic setting and use its solutions as descent directions in our augmented Lagrange algorithm for the numerical solution of our shape optimization problem.

\subsection{Optimization Problem}
\begin{figure*}[htp]
    \centering
    \begin{tikzpicture}[scale=1.0]
	\draw [thick, color=black, fill=gray!20] (0,0) rectangle (6,3);
	
	\draw [fill=white] (3,1.5) circle [radius=0.5];
	\draw [-latex] (2.64644,1.85355) -- (2.2929,2.2071) node [above] {$\mathbf{n}$};
	
	\draw (0.5,3.4) node {$B$};
	\draw (1,2) node {$\Omega$};
	\draw (2.9,1.6) node {$E$};
	\draw (-0.3,1.5) node {$\Gamma_{\mathrm{in}}$}; 
	\draw (6.55,1.5) node {$\Gamma_{\mathrm{out}}$};
	\draw (3,-0.4) node {$\Gamma_{\mathrm{slip}}$};
	\draw (3,3.4) node {$\Gamma_{\mathrm{slip}}$};
	\draw [tips, {latex}-](3.33,1.83) to [bend right] (3.7,2.1) node [above] {$\Gamma$};
	
	\draw (-1.25,3.2) node {$\mathbf{v}_{\infty}$};
	\foreach \y in {0,0.2,...,3} {
		\draw [-latex] (-1.5,\y) -- (-0.8,\y);
	};
	\draw (-1.5,0) -- (-1.5,3);
	
	\draw (6.9,0) -- (6.9,3);
	\foreach \y in {0, 0.2,...,3}{
		\pgfmathsetmacro{\x}{7.6-0.3*exp(-3.0*(\y-1.5)^2.0}
		\draw [-latex] (6.9,\y) -- (\x,\y);
	};
\end{tikzpicture}
    \caption{Schematic sketch of the computational domain.}
    \label{fig:domain}
\end{figure*}
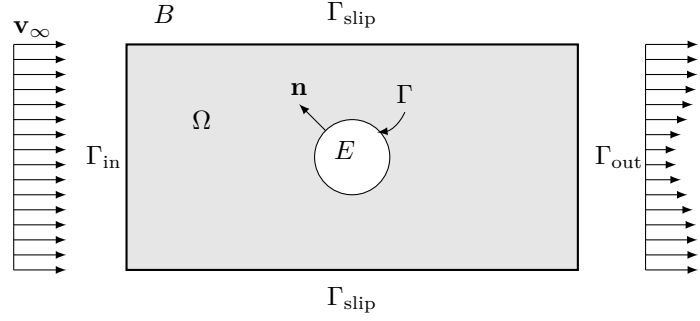

We now introduce the mathematical setting for our hydrodynamic shape optimization problem, where we refer to Fig. \ref{fig:domain} for the geometrical setting and the notation. The governing equations here are either given by the stationary Navier-Stokes equations for incompressible fluids and laminar flows or the Reynolds averaged Navier-Stokes (RANS) equations which we consider for turbulent flow at high Reynolds numbers. Note that in the later case turbulence modeling is required. Our aim is to find the shape of a generic obstacle $E \subset B$ with Lipschitz boundary located within the flow channel $B \subset \mathbb{R}^d$ which has minimal drag. The state then is given by the velocity $\mathbf{v} \colon \Omega \rightarrow \mathbb{R}^d$ and pressure $p \colon \Omega \rightarrow \mathbb{R}$, which are assumed to be unique on $\Omega$, and thus a mapping $\Omega \mapsto \mathbf{y}(\Omega) = (\mathbf{v}, p)$ exists.
It is considered that the domain $\Omega \subset B$ is completely filled with fluid which flows in from the left boundary $\Gamma_{\mathrm{in}}$ with a prescribed velocity $\mathbf{v}_\infty$ towards the outflow boundary $\Gamma_{\mathrm{out}}$ on the right. The top and bottom boundaries $\Gamma_{\mathrm{slip}}$ are considered as frictionless slip wall boundaries. 
For all test cases, the whole boundary $\Gamma$ of the obstacle $E$ is considered to be the control of the total force acting on its boundary, viz.
\begin{equation}\label{Forcefunctional}
    \mathbf{F}(\Omega) = \int_{\Gamma \subset \partial \Omega} \big( \mu (\nabla \mathbf{v} + {\nabla \mathbf{v}}^T) - \mathbf{I}p \Big) \cdot \mathbf{n} \, \mathrm{d} s,
\end{equation}
where $\mathbf{I}$ is the identity tensor. However, additionally to the flow state also geometric constraints need to be considered. On the one hand, the geometry has to be held in place to prevent the body from moving during the optimization process. On the other hand the geometry would shrink to a single point if the volume is not preserved. Thus, in the problem formulation the location of the barycenter and the wetted volume $\Omega$ are prescribed to be the vector $\bar{ \mathbf{b}}\in\mathbb{R}^d$ and the constant $\bar c>0$, respectively. Both geometrical constraints 
supplement the PDE constraints and the optimization problem reads 
\begin{subequations}
\begin{equation}
    \min_{\Omega}  \quad - \frac{\mathbf{v}_\infty}{\| \mathbf{v}_\infty \|} \cdot \mathbf{F}(\Omega) 
    \label{eq:optProblemObjectiveFunction}
\end{equation}
subjected to 
\begin{align}
    \rho (\mathbf{v} \cdot \nabla) \mathbf{v} - \div{\big(\mu (\nabla \mathbf{v} + {\nabla \mathbf{v}}^T)\big)} &= - \nabla p & \text{in} \, \Omega \label{eq:optProblemMomentumEquation}, \\
    - \nabla \cdot \mathbf{v} &= 0 & \text{in} \, \Omega, \label{eq:optProblemContinuityEquation} \\
    \mathbf{v} &= \mathbf{0} & \text{on} \, \Gamma, \label{eq:optProblemNoSlipWallBC} \\
    \mathbf{v} &= \mathbf{v}_\infty &  \, \text{on} \, \Gamma_{\mathrm{in}}, \label{eq:optProblemInletBC} \\
    \mathbf{v} \cdot \mathbf{n} = 0, \, \mathbf{n} \cdot \boldsymbol{\tau} \cdot \mathbf{t} &= \mathbf{0} & \text{on} \, \Gamma_{\mathrm{slip}}, \label{eq:optProblemSlipWallBC} \\
    \mu (\nabla \mathbf{v} + {\nabla \mathbf{v}}^T) \cdot \mathbf{n} &= p \mathbf{n} & \text{on} \, \Gamma_{\mathrm{out}}, \label{eq:optProblemOutletBC} \\
    \mathbf{b} = \frac{\int_{\Omega} \mathbf{x} \, \mathrm{d} x}{\int_{\Omega} 1 \, \mathrm{d} x}  - \bar{\mathbf{b}} &= \mathbf{0}, &  \label{eq:optProblemBarycenterConstraint} \\ 
    c = \int_{\Omega} 1 \, \mathrm{d} x - \bar{c} &= 0, & \label{eq:optProblemVolumeConstraint}
\end{align}
\end{subequations}
where the local tangent direction is defined by $\mathbf{t} = \mathbf{t}_\tau / \| \mathbf{t}_{\tau} \|_2$ with the tangent projection of the shear force $ \mathbf{t}_{\tau} \coloneqq \boldsymbol{\tau} \cdot \mathbf{n} - (\mathbf{n}^T \boldsymbol{\tau} \mathbf{n}) \mathbf{n}$ and $\boldsymbol{\tau} = \big(\mu (\nabla \mathbf{v} + {\nabla \mathbf{v}}^T) - \mathbf{I}p \big)$,  $\mathbf{v}_\infty \in \mathbb{R}^d$ denotes the inflow velocity, and $\rho, \mu \in \mathbb{R}_+$ are the density and viscosity of the fluid. The vector containing the geometrical constraints from \eqref{eq:optProblemBarycenterConstraint}-\eqref{eq:optProblemVolumeConstraint} is given by $\mathbf{g} = (\mathbf{b}, c)^T$, where $\mathbf{b}$ and $c$ refer to barycenter and volume residuals, respectively. To deal with the constraint problem above, a common approach is to introduce Lagrange multipliers. Therefore, we obtain the Lagrangian of problem \eqref{eq:optProblemObjectiveFunction}-\eqref{eq:optProblemVolumeConstraint} as
\begin{equation}
	\begin{split}
    \mathcal{L}(\Omega, \mathbf{v}, p, \hat{\mathbf{v}}, \hat{p}, \boldsymbol{\lambda}, \boldsymbol{\lambda}_b, \lambda_c) = 
    - \frac{\mathbf{v}_\infty}{\| \mathbf{v}_\infty \|} \cdot \mathbf{F}(\Omega) \\
    + \int_{\Omega} \hat{\mathbf{v}} \cdot \Big( \rho (\mathbf{v} \cdot \nabla) \mathbf{v} - \div{\big(\mu (\nabla \mathbf{v} + {\nabla \mathbf{v}}^T)\big)} + \nabla p \Big) \, \mathrm{d} x \\
    - \int_{\Omega} \hat{p} (\nabla \cdot \mathbf{v}) \, \mathrm{d} x 
    + \int_{\Gamma} \boldsymbol{\lambda} \cdot \mathbf{v} \, \mathrm{d} s 
    + \boldsymbol{\lambda}_b \cdot \mathbf{b} + \lambda_c c + \frac{\rho_b}{2} \| \mathbf{b} \|^2 + \frac{\rho_c}{2} c^2,
	\end{split}
    \label{eq:augmentedLagrangianFunction}
\end{equation}
where $\rho_{b}, \rho_c \in \mathbb{R}_{+}$ are sufficiently large penalty factors. The variables $(\hat{\mathbf{v}}, \hat{p})$ are the Lagrange multipliers for the PDE constraints \eqref{eq:optProblemMomentumEquation} and \eqref{eq:optProblemContinuityEquation} which also are associated with the adjoint state and $\boldsymbol{\lambda} \in \mathbb{R}^d$ is a Lagrange multiplier for considering the Dirichlet boundary condition \eqref{eq:optProblemNoSlipWallBC} that holds on the deformed boundary $\Gamma$. 
Note that when considering RANS equations the turbulence model implies the solution of additional state equations and thus corresponding Lagrange multipliers. Here we assume the influence of turbulence effects to be small which justifies the frequently employed frozen turbulence assumption, cf. \cite{dwight2006effects, othmer2008continuous, stuck2013adjoint}.
The multipliers $\boldsymbol{\lambda}_b \in \mathbb{R}^d$ and $\lambda_c \in \mathbb{R}$, belonging to the geometric constraints \eqref{eq:optProblemBarycenterConstraint} and \eqref{eq:optProblemVolumeConstraint}, are determined with an augmented Lagrange method and considered to be constant during the shape optimization process, as described below in detail. Following the standard approach, as described in textbooks like \cite{hinze2008, ulbrich2012}, the penalized objective function
\begin{equation}
    J(\Omega) = - \frac{\mathbf{v}_\infty}{\| \mathbf{v}_\infty \|} \cdot \mathbf{F}(\Omega) + \frac{\varrho_b}{2} \| \mathbf{b} \|^2 + \frac{\varrho_c}{2} c^2
    \label{eq:penalisedObjectiveFunction}
\end{equation}
can be expressed by the Lagrangian \eqref{eq:augmentedLagrangianFunction} and therewith the constraint problem \eqref{eq:optProblemObjectiveFunction}~-~\eqref{eq:optProblemVolumeConstraint} can be transformed into an unconstrained problem
\begin{equation}
    \min_\Omega \quad \sup_{(\hat{\mathbf{v}}, \hat{p})} \mathcal{L}(\Omega, \mathbf{v}, p, \hat{\mathbf{v}}, \hat{p}, \boldsymbol{\lambda}, \boldsymbol{\lambda}_b, \lambda_c).
    \label{eq:unconstrainedShapeOptimisationProblem}
\end{equation} 
The first order optimality system for \eqref{eq:unconstrainedShapeOptimisationProblem} is given by the partial derivative of \eqref{eq:augmentedLagrangianFunction} with respect to each argument. The \textit{Fr{\'e}chet derivative}, w.r.t. the adjoint state $(\hat{\mathbf{v}}, \hat{p})$ gives the boundary value problem \eqref{eq:optProblemMomentumEquation}~-~\eqref{eq:optProblemOutletBC}. Differentiation of \eqref{eq:augmentedLagrangianFunction} with respect to the state $(\mathbf{v}, p)$ and integration by parts leads to the adjoint equation system
\begin{subequations}
\begin{align}
    - \div{ \big( \mu (\nabla \hat{\mathbf{v}} + {\nabla \hat{\mathbf{v}}}^T) \big)} - \rho (\mathbf{v} \cdot \nabla) \hat{\mathbf{v}} + \rho {\nabla \mathbf{v}}^T \, \hat{\mathbf{v}} &= - \nabla \hat{p} & \text{in} \, \Omega, \label{eq:adjointMomentumEquation} \\
    \nabla \cdot \hat{\mathbf{v}} &= 0 & \text{in} \, \Omega, \label{eq:adjointContinuityEquationBC} \\
    \hat{\mathbf{v}} + \frac{\mathbf{v}_\infty}{\| \mathbf{v}_\infty \|} &= \mathbf{0} & \text{on} \, \Gamma, \label{eq:adjointDesignShapeBC}\\
    \hat{\mathbf{v}} &= \mathbf{0} & \text{on} \, \Gamma_{\mathrm{in}} \label{eq:adjointNonDesignBC} \\
   \mu (\nabla \hat{\mathbf{v}} + {\nabla \hat{\mathbf{v}}}^T) \cdot \mathbf{n} &= \hat{p} \mathbf{n} - \rho (\mathbf{v}\cdot\mathbf{n})\hat{\mathbf{v}} & \text{on} \, \Gamma_{\mathrm{out}} \label{eq:adjointOutletBC} \\
    \hat{\mathbf{v}} \cdot \mathbf{n} = 0, \, \mathbf{n} \cdot \hat{\boldsymbol\tau}\cdot \mathbf{t} &= \mathbf{0} & \text{on} \, \Gamma_{\mathrm{slip}} \label{eq:adjointSlipwallBC}
\end{align}
and with
\begin{equation}
    \boldsymbol{\lambda} = \hat{p} \mathbf{n} - \mu (\nabla \hat{\mathbf{v}} + {\nabla \hat{\mathbf{v}}}^T) \cdot \mathbf{n}
\end{equation}
\end{subequations}
where $\hat{\boldsymbol\tau} = \mu (\nabla \hat{\mathbf{v}} + {\nabla \hat{\mathbf{v}}}^T) - \mathbf{I}\hat{p}$.
Together with the PDE constraints of the optimization problem \eqref{eq:optProblemMomentumEquation}~-~\eqref{eq:optProblemOutletBC} the adjoint equation system~\eqref{eq:adjointMomentumEquation}~-~\eqref{eq:adjointSlipwallBC} characterizes a saddle point of the Lagrangian \eqref{eq:augmentedLagrangianFunction}, which is assumed to be a unique stationary point of $\mathcal{L}$. In order to compute the shape derivative of  the augmented Lagrangian \eqref{eq:augmentedLagrangianFunction} the domain $\Omega$ is made variable by a family of transformations $\{\mathbf{T}_t\}_{t \geq 0}$ with the parameterized perturbation of identity
\begin{equation}
	\begin{split}
    \mathbf{T}_t = \id + t \mathbf{u}_p \\
    \text{with} \quad \mathbf{T}_0(\Omega) = \Omega \; \text{and} \quad \mathbf{u}_p \in V^{1,p}_0(\Omega),
	\end{split}
    \label{eq:domainTransformation}
\end{equation}
where $V_0^{1,p} := \lbrace \mathbf{u} \in W^{1,p}(\Omega, \mathbb{R}^d) : \mathbf{u} = \mathbf{0} \text{ a.e. on } \partial \Omega \setminus \Gamma \rbrace$, and $t \geq 0$ is a step size, cf. \cite{kuhl2019decoupling}. Application of the first order optimality condition
\begin{equation}
    \frac{\partial}{\partial t} \mathcal{L}(\mathbf{T}_t(\Omega), \mathbf{v}, p, \hat{\mathbf{v}}, \hat{p}, \boldsymbol{\lambda}, \boldsymbol{\lambda}_b, \lambda_c) \Big|_{t = 0} = 0 \label{eq:firstOrderOptimalityCondition}
\end{equation}
leads to the surface formulation of the shape derivative of the objective function $J$, compare e.g. \cite{BZ93}. For its representation we define
\begin{equation}
		\gamma \coloneqq - \mu \frac{\partial \hat{\mathbf{v}}}{\partial n} \cdot \frac{\partial \mathbf{v}}{\partial n}
		+  \boldsymbol{\lambda}_b \cdot \frac{\mathbf{x} - \boldsymbol{\beta}}{\int_{\Omega} 1 \, \mathrm{d} x} + \lambda_c + \rho_{b} \frac{(\boldsymbol{\beta} - \bar{\mathbf{b}})(\mathbf{x} - \boldsymbol{\beta})}{\int_{\Omega} 1 \, \mathrm{d} x} 
		+ \rho_{c} \big( \int_{\Omega} 1 \, \mathrm{d} x - \bar{c} \big)
\end{equation}
where
\begin{equation}
	\boldsymbol{\beta} = \frac{\int_{\Omega} \mathbf{x} \, \mathrm{d} x}{\int_{\Omega} 1 \, \mathrm{d} x}
\end{equation}
is the barycenter of $\Omega$, c.f. \cite{bello1997differentiability, schulz2016computational}. Then
\begin{gather}
	J^\prime (\Omega; \mathbf{u}) = 
	\int_{\Gamma}
	\gamma \mathbf{u} \cdot \mathbf{n}
	\, \mathrm{d} s, \label{eq:optProblemShapeDerivative}
\end{gather}
is the shape derivative of $J$ at $\Omega$ into direction of $\mathbf{u}$. While the Lagrange multipliers $(\hat{\mathbf{v}}, \hat{p})$ are given by the solution of the boundary value problem \eqref{eq:adjointMomentumEquation}~-~\eqref{eq:adjointSlipwallBC}, the multipliers $\boldsymbol{\lambda}_b$ and $\lambda_c$ for barycenter and volume constraint are not given by the solution of a system of PDEs and rather have to be determined by an augmented Lagrange method like outlined in algorithm~\ref{alg::augmented_lagrange} of the next section, c.f. \cite{andreani2008augmented, allaire2020}.
Hence, with the shape derivative from \eqref{eq:optProblemShapeDerivative} the minimization problem in \eqref{eq:pLaplacianVariationalForm} with $W^{1,p}(\Omega,\mathbb{R}^d)$ replaced by $V_0^{1,p}(\Omega,\mathbb{R}^d)$ is associated with the variational form to the $p$-Laplacian problem
\begin{equation}
    \left.\begin{array}{ll}
    - \div{\big( (\nabla \mathbf{u}_p \colon \nabla \mathbf{u}_p)^{\frac{p-2}{2}} \nabla \mathbf{u}_p\big)} = 0 & \text{ in } \, \Omega, \\
    (\nabla \mathbf{u}_p \colon \nabla \mathbf{u}_p)^{\frac{p-2}{2}} \frac{\partial \mathbf{u}_p}{\partial n} = \gamma \mathbf{n} & \text{ on } \, \Gamma, \\
    \mathbf{u}_p = 0 & \text{ on } \, \partial \Omega\setminus\Gamma.
    \end{array}\right\}
    \label{eq:plaplacianProblem}
\end{equation}
Note that $\Gamma \subset \partial \Omega$ is the part of the boundary which is free for deformation. 
This approach now can also be interpreted as a generalization by the $p$-Laplace setting of the approach proposed in \cite{schulz2016computational} where the descent direction is defined regarding a Steklov-Poincar{\'e} operator; in other words the Dirichlet-to-Neumann map is applied to $\gamma \mathbf{n}$. Therewith, the descent direction in \cite{allaire2004structural, schulz2016computational} is obtained by solving an elliptic boundary value problem similar to
\begin{equation}
    \left.\begin{array}{ll}
       - \Delta \mathbf{u} = 0 & \text{ in } \, \Omega \\
        \frac{\partial \mathbf{u}}{\partial n} = \gamma \mathbf{n} & \text{ on } \, \Gamma \\
        \mathbf{u} = 0 & \text{ on } \partial \Omega \backslash \Gamma
    \end{array}\right\}
    \label{eq:LaplaceProblem}
\end{equation}
which is the Euler-Lagrange equation of the minimization problem
\begin{equation}
    \min_{\{\mathbf{u} \in H^1(\Omega, \mathbb{R}^d); \mathbf{u} = 0 \text{ on } \partial\Omega\setminus \Gamma\}} \quad \frac{1}{2} \int_{\Omega} \nabla \mathbf{u} \colon \nabla \mathbf{u} \, \mathrm{d} x
    - \int_{\Gamma} \gamma \mathbf{n} \cdot \mathbf{u} \, \mathrm{d} s,
\end{equation}
This obviously represents a special case of \eqref{eq:pLaplacianVariationalForm} with $p = 2$ and serves as a reference for the numerical experiments discussed in section~\ref{sc:numerics}. Let us also note that this concept was enhanced in \cite{onyshkevych2020mesh} where a nonlinear extension operator has been introduced. The observation of strong distortion of the discrete grid within the main deformation direction motivated adding a nonlinear advection term to the PDE in \eqref{eq:LaplaceProblem}. Therewith, greater deformations are possible and the mesh quality is reasonable even after large deformations of the initial grid.
However, as our numerical results show, the $p-$Laplace relaxation of problem  \eqref{eq:infinityLaplacianProblem} provides a systematic approach to fluid dynamic shape optimization also guaranteeing meshes of high quality after large deformations, in particular when the optimal shape has sharp corners.

\section{Optimization Algorithm} \label{sc:algorithm}
In Algorithm \ref{alg::augmented_lagrange} we specify our minimization algorithm.
\begin{algorithm}[!ht]
	\caption{Augmented Lagrange Optimization}
    \begin{algorithmic}[1]
		\Require{$\Omega^0 \subset B$, $\bar{\mathbf{b}}$, $\boldsymbol{\lambda}_{b} \in \mathbb{R}^d$, $\lambda_c \in \mathbb{R}$, $\bar{c}$, $\varrho_{b}$, $\varrho_{c}$, $\varrho^{\mathrm{inc}}$, $\tau_{b}$, $\tau_c >0$}
		\State{$k \gets 0$}
		\Repeat
    		\State $\Omega^{k+1} \gets \underset{\Omega}{\mathrm{arg \, min}} \; \mathcal{L}(\Omega, \mathbf{v}, p, \hat{\mathbf{v}}, \hat{p}, \boldsymbol{\lambda}_{b}, \lambda_c)$ 	 \label{eq:optimisationLoop} \label{alg::augmented_lagrange::shapeOptimisation}
    		\State $\mathbf{b} \gets \int_{\Omega^{k}} \mathbf{x} \, \mathrm{d} x / \int_{\Omega} 1 \, \mathrm{d} x - \bar{\mathbf{b}}$
    		\State $c \gets \int_{\Omega^{k}} 1 \, \mathrm{d} x - \bar{c}$
    		
    		\If{$\| \mathbf{b} \| > \tau_{b}$}
    		    \State $\varrho_{b} \gets \varrho^{\mathrm{inc}} \varrho_{b}$
    		\Else
    		    \State $\boldsymbol{\lambda}_{b} \gets \boldsymbol{\lambda}_{b} + \varrho_{b} b$ \label{alg::augmented_lagrange::barycentreMultiplierUpdate}
    		\EndIf
    
    		\If{$|c| > \tau_c$}
    		    \State $\varrho_c \gets \varrho^{\mathrm{inc}}\varrho_c$
    		\Else 
    		    \State $\lambda_c \gets \lambda_c + \varrho_c c$ \label{alg::augmented_lagrange::volumeMultiplierUpdate}
    		\EndIf
    		\State $k \gets k + 1$
		\Until{$\vert J(\Omega^{k+1}) - J(\Omega^k) \vert < \epsilon \vert J(\Omega^0)\vert$}
	\end{algorithmic}
    \label{alg::augmented_lagrange}
\end{algorithm}
It was used before, in e.g. \cite{schulz2016computational} to determine the multipliers $\boldsymbol{\lambda}_b$ and $\lambda_c$ for the barycenter and volume constraints of a shape optimization problem. Other than in \cite{schulz2016computational} the update for the penalty factors $\rho_b, \rho_c \in \mathbb{R}_+$ and $\boldsymbol{\lambda}_b \in \mathbb{R}^d$ and $\lambda_c \in \mathbb{R}$, respectively, is separated which is motivated by the big difference in the order of magnitude between both penalty factors ($\varrho_b \propto 10^8$ and $\varrho_c \propto 10^2$) for the numerical experiments in this paper. However, the critical part of algorithm~\ref{alg::augmented_lagrange} refers to step~\ref{eq:optimisationLoop} which is therefore farther explored in algorithm~\ref{alg::shape}. To compute the multipliers $\boldsymbol{\lambda}_b$ and $\lambda_c$ the shape optimization algorithm~\ref{alg::shape} is called in a sequence of convergence tolerances $\epsilon \in \mathbb{R}_+$, c.f. step~\ref{alg::shape::convergenceCriteria} in algorithm~\ref{alg::shape}, in order to keep the computation numerically stable and circumvent unfeasible shapes during the optimization. 
\begin{algorithm}[!ht]
	\caption{Shape Optimization}
    \begin{algorithmic}[1]
        \Require{$\Omega$, $\bar{\mathbf{b}}$, $\boldsymbol{\lambda}_b$, $\bar{c}$, $\rho_c$, $\lambda_c$, $\epsilon$}
    	\State{$i \gets 0$}
    	\Repeat
       		\State solve primal equations \eqref{eq:optProblemMomentumEquation}~-~\eqref{eq:optProblemOutletBC} for $\mathbf{v},p$ \label{alg::shape::solvePrimal}
       		\State solve adjoint equations \eqref{eq:adjointMomentumEquation}~-~\eqref{eq:adjointSlipwallBC} for $\hat{\mathbf{v}}, \hat{p}$ \label{alg::shape::solveAdjoint}
       		\State compute $\gamma$ according to \eqref{eq:optProblemShapeDerivative} \label{alg::shape::computeShapeSensitivity}
       		\State solve $p$-Laplace problem \eqref{eq:plaplacianProblem} \label{alg::shape::solvePLaplace}
       		\State move grid points according to \eqref{eq:domainTransformation} \label{alg::shape::moveGrid}
       		\State $i \gets i + 1$
    	\Until{$\vert J(\Omega^{i+1}) - J(\Omega^i) \vert < \epsilon \vert J(\Omega^0)\vert$} \label{alg::shape::convergenceCriteria}
	\end{algorithmic}
    \label{alg::shape}
\end{algorithm}

A conventional, pressure-based, second-order accurate finite-volume scheme for a cell-centered variable arrangement is employed to discretize the
partial differential equations of the primal \eqref{eq:optProblemMomentumEquation}~-~\eqref{eq:optProblemOutletBC} and adjoint systems \eqref{eq:adjointMomentumEquation}~-~\eqref{eq:adjointSlipwallBC}, cf. \cite{rung2009challenges, stuck2013adjoint, kuhl2020adjoint}. The existing infrastructure and generic subroutines of the fluid solver can be re-used with limited effort for the implementation. In case of $p = 2$, the initial guess with $\mathbf{u}_p = 0$ in $\Omega$ leads to convergence of the implementation. For $p > 2$ a sequence of problems in $p$ is solved in order to obtain an initial guess for the discrete problem with the desired value for $p$. Within the procedure the solution of the preceding smaller $p$ is used as initial guess. Approximations that offer an initial guess only experience a smaller convergence tolerance in order to save computational time.

\section{Applications}
\label{sc:numerics}
Three fluid dynamic  applications are discussed to investigate the performance of the $p$-Laplace approach  \eqref{eq:pLaplacianVariationalForm} in shape optimization. All cases are concerned with drag minimization at steady state and subjected to conserve the wetted volume and its barycenter, c.f. Fig. \ref{fig:domain}. Starting from the Laplace expression with $p=2$, the first case analyses the influence of increasing $p$ for a frequently referenced $2D$ Stokes flow example that features a pointed oval optimum. Emphasis is given to \begin{inparaenum}[(a)] \item the convergence of the optimization, \item the final shape and the attainable drag reduction as well as \item the quality of the mesh updates\end{inparaenum}. The second case demonstrates the applicability of the approach for an analogue $3D$ configuration and analyses the same aspects (a-c). Since both initial applications refer to low Reynolds (low Re) number flows, a third example is supplemented to scrutinize the performance in a $2D$ turbulent flow at high Reynolds number. 

\subsection{Drag Optimization in $2D$ Low Re Flow} 
\label{sc:application::2DLowRe}
The first case studies the drag minimization of a $2D$ circular cylinder exposed to low Reynolds number flow for $p = 2, 3, \ldots, 6$. The computational domain is illustrated in Fig. \ref{fig:domain}. The case is related to the setting initially described by Pironneau \cite{pironneau1973optimum}. Instead of the Stokes flow considered in \cite{pironneau1973optimum}, we employ a low Re Navier-Stokes formulation of the boundary value problem. 

The initial cylinder features a unit diameter ($D=1 [m]$) and is centered in a channel of length $50 [m]$ and height $10 [m]$. The fluid is characterized by a unit density and dynamic viscosity, i.e.  $\rho = 1 [kg/m^3]$, $\mu = 1 [Ps \cdot s]$. Dirichlet conditions are imposed at the inflow $\Gamma_{\mathrm{in}}$ with $\mathbf{v}_\infty = (1, 0)^T [m/s]$, which yields a unit Reynolds number $\mathit{Re}_D = 1$. Slip wall boundary conditions are applied to the top and bottom boundaries of the channel. Outflow boundary conditions \eqref{eq:optProblemOutletBC} are used along  $\Gamma_{\mathrm{out}}$.
\begin{figure}[!ht]
    \centering
    \includegraphics[width=0.7\textwidth, trim=400 200 400 200, clip]{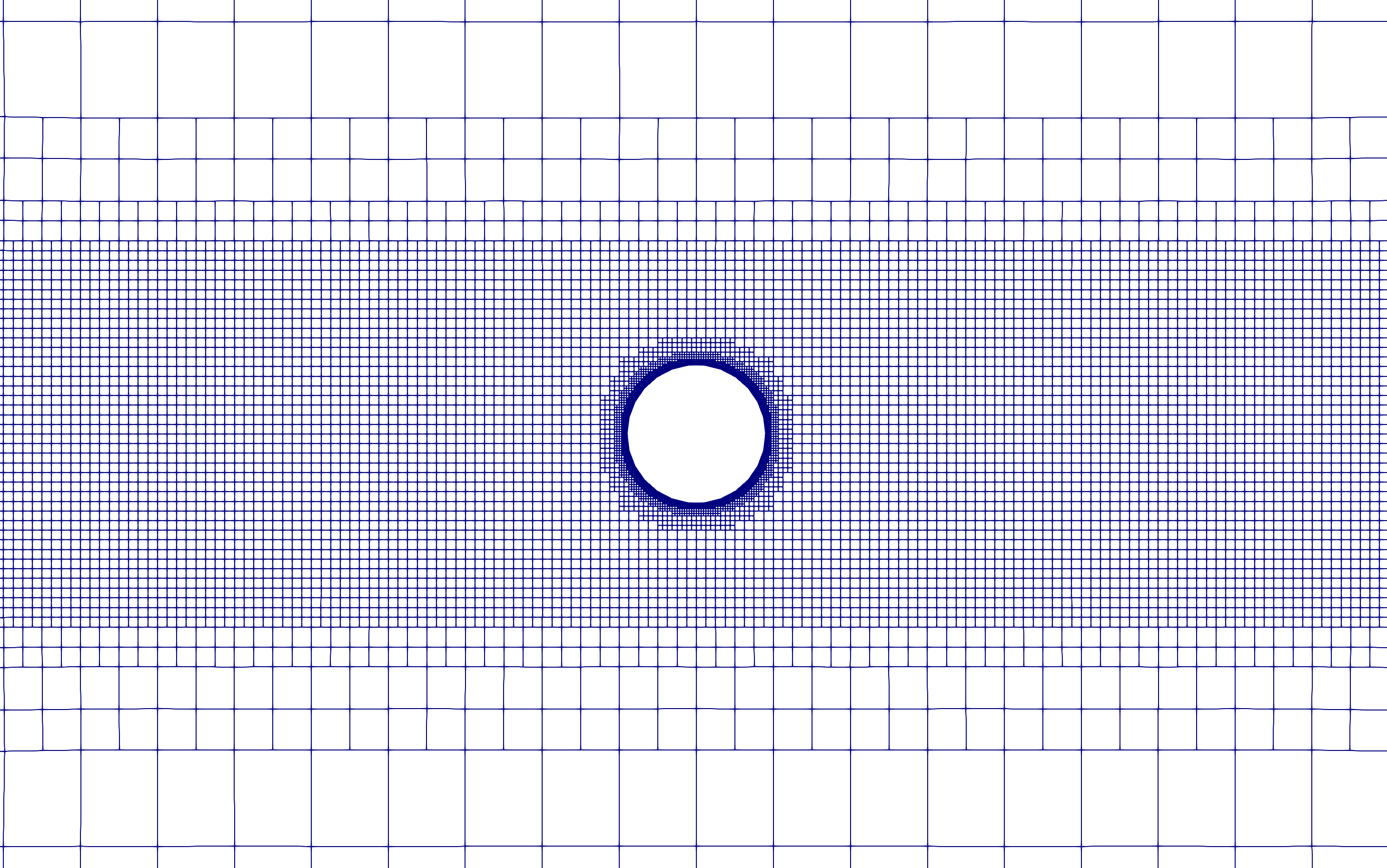}
    \caption{Initial mesh of the $2D$ low Reynolds number problem.}
    \label{fig:cylinder2DInitialMesh}
\end{figure}
The initial grid is displayed in Fig. \ref{fig:cylinder2DInitialMesh} and features $927$ evenly distributed cells along the circumference of the cylinder boundary $\Gamma$. The CAD-free optimization procedure employs an unstructured, locally refined mesh with approximately $14.5$k control volumes which is deformed along with the surface using \eqref{eq:domainTransformation}. Mind that the  optimal shape is expected to reveal pointy tips at the front and the aft, which shall deliberately not be preempted by the initial discretization of the boundary $\Gamma$. The target barycentre and wetted volume are set to $\bar{\mathbf{b}} = (0, 0)^T$ and $\bar{c} = 50 \cdot 10 - \pi / 4 [m^2]$ respectively. All parameters employed to initialize the algorithms \ref{alg::augmented_lagrange} and \ref{alg::shape} are displayed in Table \ref{tab:augmentedLagrangeInput}. The multipliers for the barycenter and volume are initialized with $\boldsymbol{\lambda}_b = (0, 0)^T$ and $\lambda_c = 0$, respectively. The sequence of tolerances applied to the convergence criteria of the shape optimisation problem is given    by $\epsilon = 10^{-1}, 10^{-2}, \ldots, 10^{-7}$. However, it turned out that the multipliers $\boldsymbol{\lambda}_b$ and $\lambda_c$ converge very fast and after four augmented Lagrange steps the multipliers are determined sufficiently exact in each computation with different values for $p$. The field values for the fixed point iteration for $p = 2$ can be initialized with $\mathbf{u}_p = 0$ for all discrete points. Larger $p$-values were initialized by solutions obtained from the previous smaller value, each to a suitable tolerance to provide an initial guess for the next following $p$-Laplace problem. Although the theory suggests to drive $p \to \infty$, lower $p$-values are of interest for large-scale applications due to the more exhausting computational effort. The numerical effort for solving the $p$-Laplace problem is investigated in \cite{loisel2020} for different values for $p$ which show that it is of polynomial complexity but depends on $p$ and the number of unknowns. The related experience from this investigation reveals an increase of computing time $T_{p}$ for one iteration of the shape optimization algorithm~\ref{alg::shape} by $T_{p = 4} / T_{p = 3} \approx 2.8$ and $T_{p = 5} / T_{p = 4} \approx 1.6$. In addition the representable floating point arithmetic of the machine limits the value for $p$.
\begin{table}[!ht]
    \centering
    \begin{tabular}{c|c c c}
        Parameter                   & $2D$ low $\mathit{Re}$    & $3D$ low $\mathit{Re}$   &  $2D$ high $\mathit{Re}$   \\ 
        \hline
        $\varrho_b$                 & $5 \cdot 10^7$            & $1 \cdot 10^3$           & $1.2 \cdot 10^8$             \\
        $\varrho_c$                 & $1 \cdot 10^2$            & $50$                     & $4 \cdot 10^2$             \\
        $\varrho^{\mathrm{inc}}$    & $2$                       & $2$                      & $1.2$                      \\
        $\tau_b$                    & $1 \cdot 10^{-6}$         & $1 \cdot 10^{-1}$        & $5 \cdot 10^{-5}$          \\
        $\tau_c$                    & $2 \cdot 10^{-2}$         & $1 \cdot 10^{-4}$        & $2 \cdot 10^{-2}$          \\
        $t$                         & $2 \cdot 10^{-3}$         & $5 \cdot 10^{-3}$        & $1 \cdot 10^{-3}$       
    \end{tabular}
    \caption{Initial values for the parameters of the augmented Lagrange procedure.}
    \label{tab:augmentedLagrangeInput}
\end{table}
\begin{figure*}[t]
    \centering
    \begin{tikzpicture}
    \begin{axis}[
        height = 0.3\textheight,
        width = 0.8\textwidth,
        xmin = 0, xmax = 300,
        ymin = 0.9, ymax = 1.0,
        grid = major,
        cycle list name = my color list,
        legend pos = outer north east,
        ylabel = $J(\Omega) / J(\Omega^0)$,
        xlabel = optimisation steps,
        legend pos = north east,
        every axis plot post/.append style={mark=none, line width=1.0pt},
    ]
    
    \addplot+ [] table [x expr=\thisrow{i}-810, y=Jn] {tikz/data/cylinder_2D/objective_p02.0.dat};
    
    \addplot+ [] table [x expr=\thisrow{i}-818, y=Jn] {tikz/data/cylinder_2D/objective_p03.0.dat};
    
    \addplot+ [] table [x expr=\thisrow{i}-781, y=Jn] {tikz/data/cylinder_2D/objective_p04.0.dat};
    
    \addplot+ [] table [x expr=\thisrow{i}-355, y=Jn] {tikz/data/cylinder_2D/objective_p05.0.dat};
    
    \addplot+ [] table [x expr=\thisrow{i}-293, y=Jn] {tikz/data/cylinder_2D/objective_p06.0.dat};
    
    \legend{$p = 2$, $p = 3$, $p = 4$, $p = 5$, $p  =6$};
    
    \end{axis}
\end{tikzpicture}
    \caption{Influence of $p$-value on the evolution of the normalized objective functional.}
    \label{fig:cylinder2DobjectivePolot}
\end{figure*}
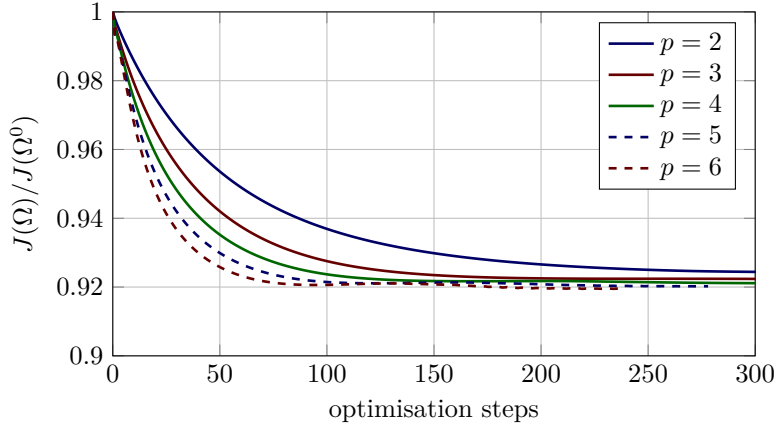

\subsubsection{Optimal Shapes and Convergence}
\begin{figure*}[t]
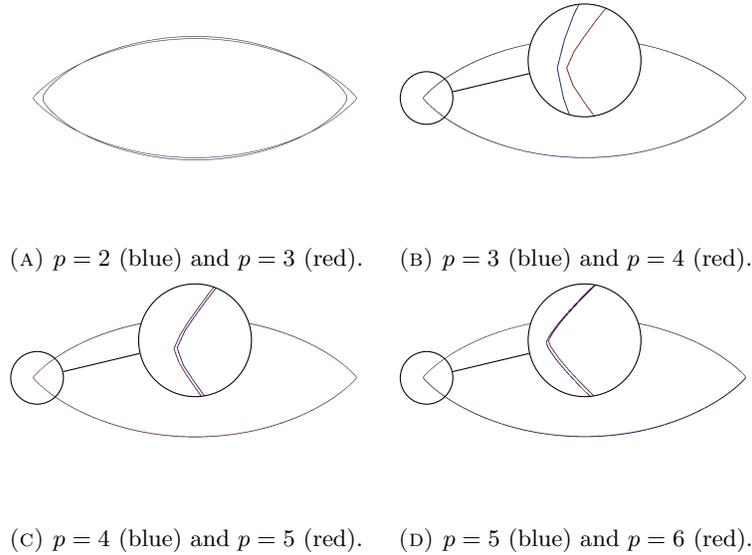

    \centering
    \begin{subfigure}[c]{0.4\textwidth}
        \input{tikz/cylinder_2D_compare_p02.0_p03.0.tikz}
        \subcaption{$p = 2$ (blue) and $p = 3$ (red).}
        \label{fig:compareFinalShapeP02P03Contours}
    \end{subfigure}
     \begin{subfigure}[c]{0.4\textwidth}
        \input{tikz/cylinder_2D_compare_p03.0_p04.0.tikz}
        \subcaption{$p = 3$ (blue) and $p = 4$ (red).}
        \label{fig:compareFinalShapeP03P04Contours}
    \end{subfigure}
    \begin{subfigure}[c]{0.4\textwidth}
        \input{tikz/cylinder_2D_compare_p04.0_p05.0.tikz}
        \subcaption{$p = 4$ (blue) and $p = 5$ (red).}
        \label{fig:compareFinalShapeP04P05Contours}
    \end{subfigure}
     \begin{subfigure}[c]{0.4\textwidth}
        \input{tikz/cylinder_2D_compare_p05.0_p06.0.tikz}
        \subcaption{$p = 5$ (blue) and $p = 6$ (red).}
        \label{fig:compareFinalShapeP05P06Contours}
    \end{subfigure}
    \caption{Influence of the $p$-value on the predicted optimal (final) shapes.}
    \label{fig:cylinder2DfinalShapeContours}
\end{figure*}

Figure \ref{fig:cylinder2DobjectivePolot} shows the normalized evolution of the drag objective \eqref{eq:penalisedObjectiveFunction} over a selected number of optimization steps for the five investigated values of $p$. The baseline solution refers to $p = 2$ which is also related to the investigations in \cite{schulz2016computational}.
\begin{table}[!ht]
    \centering
    \begin{tabular}{c|c c c c}
        $p$   & Tip angle $[^\circ]$ & $a/b _{(it=50)}$ & $J/J^0$ & conv. step \\
        \hline
        $2$ &  $163.8$ & $1.4$ & $0.9243$ & $/$   \\
        $3$ &  $152.4$ & $1.6$ & $0.9222$ & $415$ \\
        $4$ &  $140.4$ & $1.8$ & $0.9211$ & $337$ \\
        $5$ &  $129.4$ & $2.0$ & $0.9202$ & $278$ \\
        $6$ &  $116.4$ & $2.1$ & $0.9195$ & $239$ \\
    \end{tabular}
    \caption{Performance indicators obtained with different $p$-values for 2D low Re case.}
    \label{tab:cylinder2DoptimisationStats}
\end{table}
Table ~\ref{tab:cylinder2DoptimisationStats} outlines a comparison of performance indicators obtained for the five investigated values of $p$. The last column refers to the maximum number of design steps that are needed to reach a sufficiently converging objective function (cf. algorithm~\ref{alg::shape}). The table reveals that the convergence improves and fewer optimization steps are needed for larger values of $p$. For $p = 2$  the optimization could not reach convergence, but terminated after $356$ steps due to grid quality issues, which are discussed below. At this step, the convergence criterion was about two orders of magnitude above the threshold. However, from a practical point of view the convergence criteria employed in this study might appear rather strict and practical applications would also reach sufficient optima for $p=2$. Table~\ref{tab:cylinder2DoptimisationStats} also displays the final objective function values, which again reveal improvements for increasing $p$-values. Drag reductions refer to about $7,6 \%$ for the Laplacian approach with $p = 2$ and increase to  approximately $8,1 \%$ for $p=6$. Improvements seen for the objective function are attributed to the more extreme deformations obtained from the $p$-Laplacian problem \eqref{eq:pLaplacianVariationalForm} with $p > 2$. In order to judge the final shape, the opening angle at the upstream tip may also serve as a measure. The interior opening angles listed in Tab.~\ref{tab:cylinder2DoptimisationStats} decrease with greater values of $p$. Hence, increasing $p$ clearly yields more pointy tips as also indicated by the  comparison of tip shapes in Fig.~\ref{fig:cylinder2DfinalShapeContours}. The Stokes flow problem investigated in \cite{pironneau1973optimum} reported an opening angle of $120^\circ$. The present results rapidly approach the reported opening angle from above.  However, for $p = 6$ the opening angle falls below the  reference value, which is attributed to the use of Navier-Stokes instead of the Stokes flow model \cite{pironneau1973optimum}. When attention is directed to the convergence speed, the half axis ratio at the $50$th design step, i.e. $a/b _{(it=50)}$, mentioned in Tab.~\ref{tab:cylinder2DoptimisationStats} may be considered as a measure to assess the convergence speed. As all simulations start with $a/b=1$, the tabulated data renders the influence of $p$-values on the ability of the $p$-Laplace approach to rapidly adjust the shape. It is also observed, that large deformations take place at an early stage of the optimization. For example, the final half axis ratio for $p=4$ refers to approximately $2.7$ and already $2/3$ of this ratio are reached after $50$ design steps. A closer inspection of Fig.~\ref{fig:cylinder2DobjectivePolot} reveals,  that the value of the objective function slightly increases after a quick descent for  $p = 5$  (step $120$-$160$) and $p = 6$ (step $90$-$140$). This is possible due to the absence of a step size control and occurs because the multipliers that control the barycenter and the displacement strictly speaking only hold for a preceding iteration. Moreover, the $p$-Laplace problem  \eqref{eq:pLaplacianVariationalForm} is not solved exactly in every iteration of algorithm~\ref{alg::shape} and preceding results are used as the initial guess for a subsequent optimization step. A similar phenomenon can be seen in results reported by Allaire et al. \cite{allaire2004structural}, who did use a similar augmented Lagrange procedure.
The contours of the final shapes are depicted in Fig.~\ref{fig:cylinder2DfinalShapeContours}, where  optimal shapes for two consecutive values of $p$ are compared with each other, i.e. for $p = 2$ with $p = 3$, $p = 3$ with $p = 4$ and so forth. The shape contours in Fig.~\ref{fig:compareFinalShapeP02P03Contours} show significant differences, in particular at the tips of the resulting geometry. The  computation with $p = 2$ clearly leads to a more rounded shape and a larger vertical extent than the other investigated $p$ values. Shapes returned by $3 \le p \le 6$ are displayed Fig.~\ref{fig:compareFinalShapeP03P04Contours}~-~\ref{fig:compareFinalShapeP05P06Contours}. Remarkably, a general difference between the respective contours is hard to identify for $p\ge3$. Thus, close-ups are used to assess the tip region. While a rounded tip region is still observed for $p = 3$, the tip becomes more pointy for $p \ge 4$. Figures ~\ref{fig:compareFinalShapeP04P05Contours} and \ref{fig:compareFinalShapeP05P06Contours} only display small differences between the  shapes obtained with $p=4,5,6$. Thus, one can assume that the predicted optimal shapes tend to be converged for $p\ge 4$.
\begin{figure*}[t]
    \centering
    \begin{subfigure}[c]{0.18\textwidth}
        \includegraphics[width=1.0\textwidth, trim=700 100 400 100, clip]{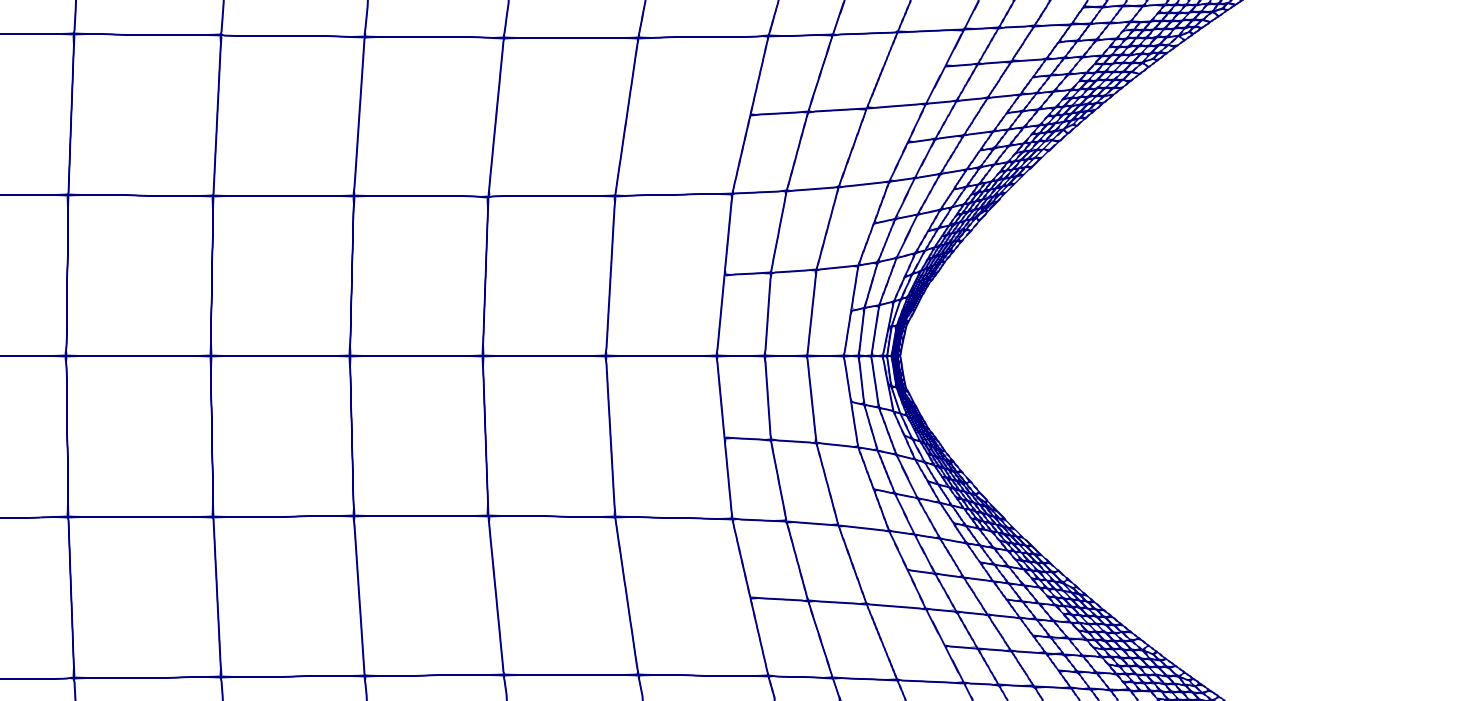}
        \subcaption{$p = 2$}
        \label{fig:cylinder2DMeshTipDetailP02}
    \end{subfigure}
    \begin{subfigure}[c]{0.18\textwidth}
        \includegraphics[width=1.0\textwidth, trim=600 100 500 100, clip]{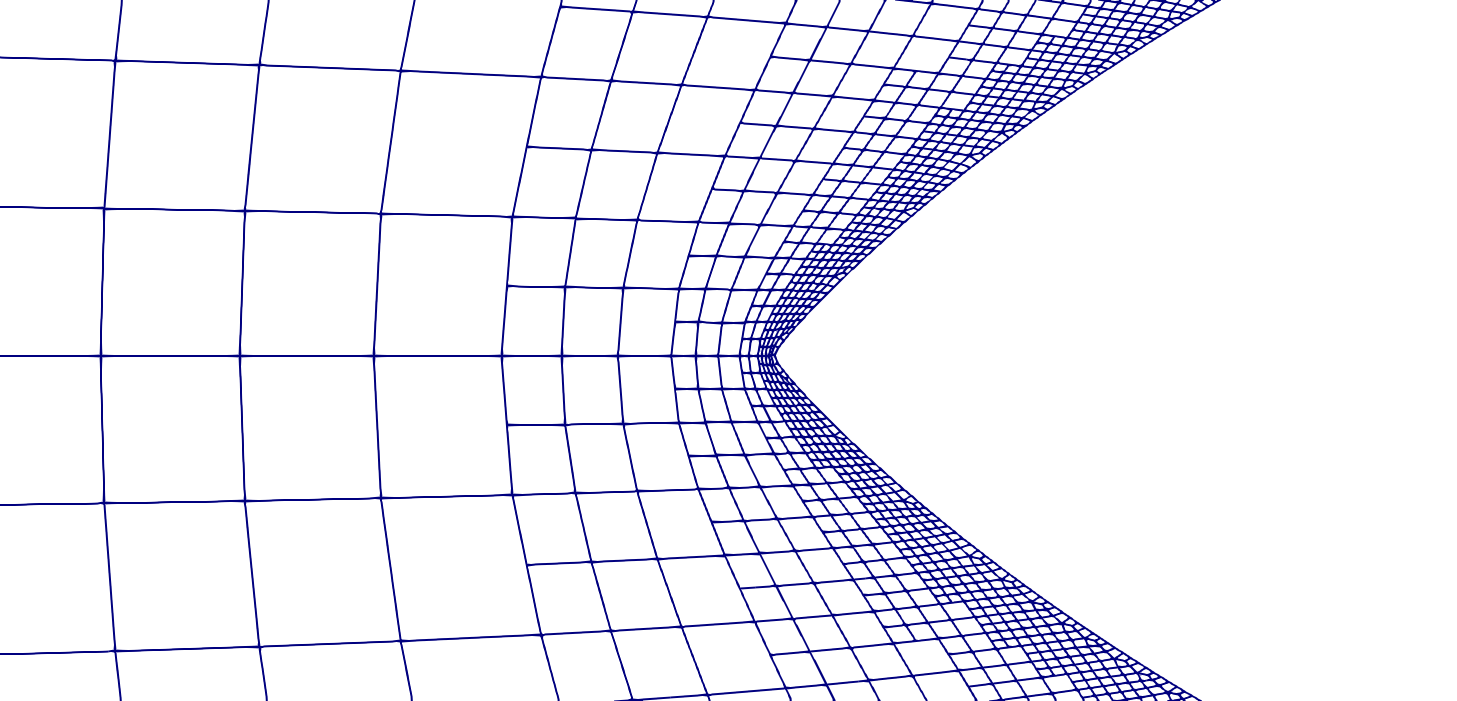}
        \subcaption{$p = 3$}
        \label{fig:cylinder2DMeshTipDetailP03}
    \end{subfigure}
    \begin{subfigure}[c]{0.18\textwidth}
        \includegraphics[width=1.0\textwidth, trim=600 100 500 100, clip]{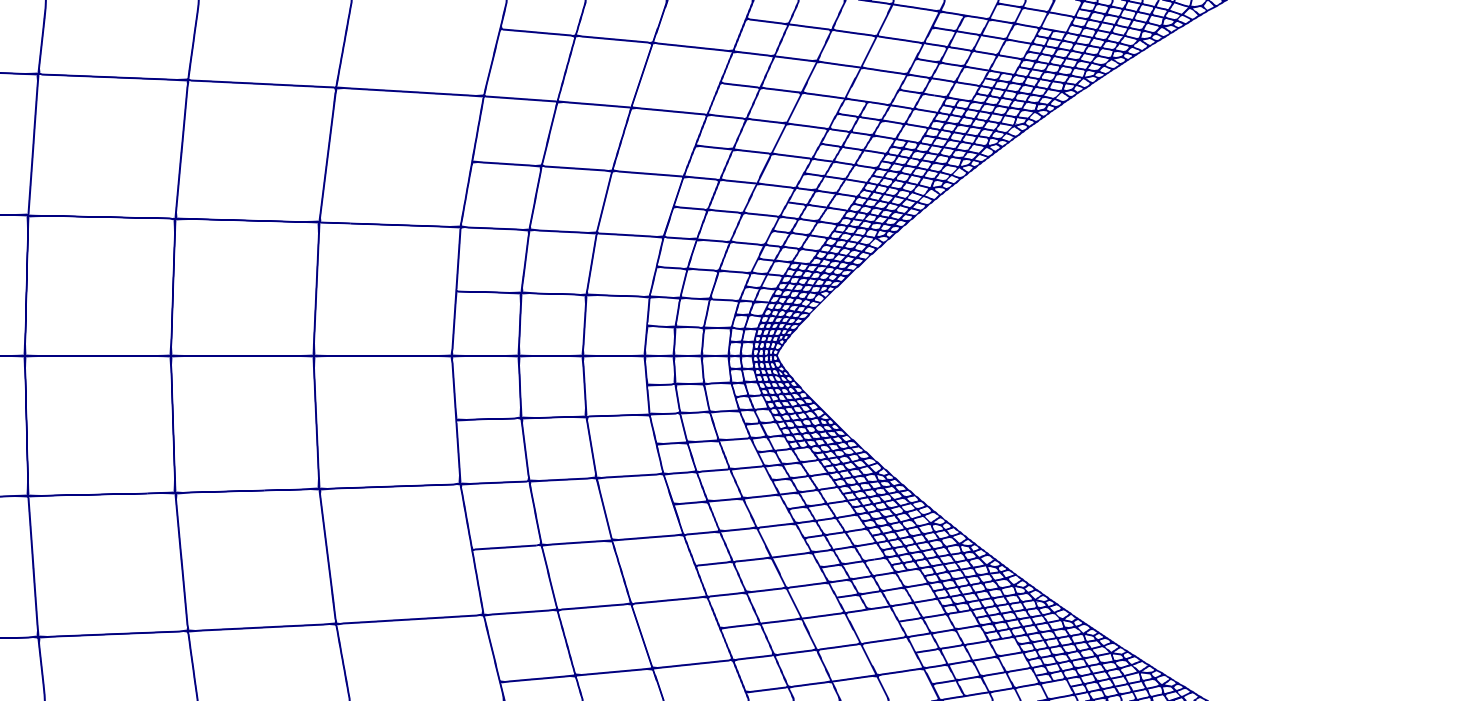}
        \subcaption{$p = 4$}
        \label{fig:cylinder2DMeshTipDetailP04}
    \end{subfigure}
    \begin{subfigure}[c]{0.18\textwidth}
        \includegraphics[width=1.0\textwidth, trim=600 100 500 100, clip]{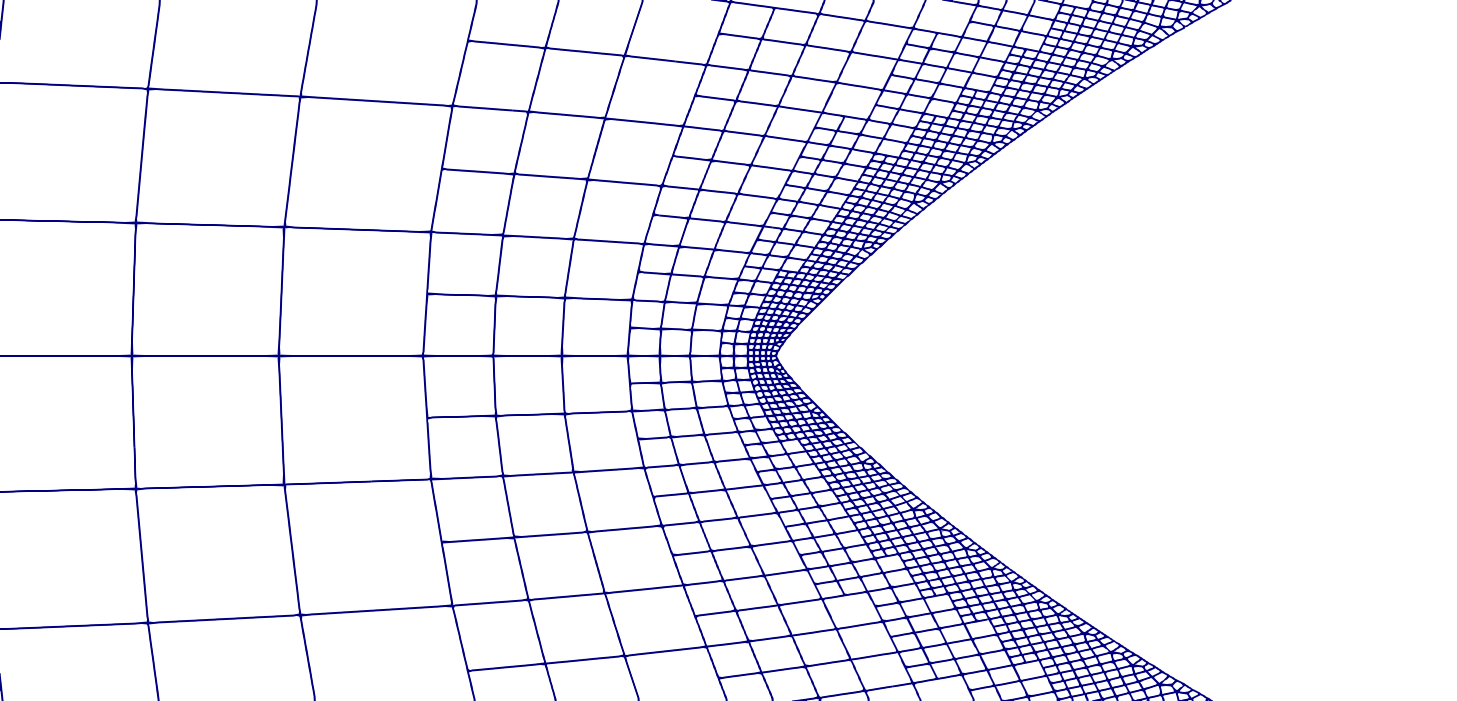}
        \subcaption{$p = 5$}
        \label{fig:cylinder2DMeshTipDetailP05}
    \end{subfigure}
    \begin{subfigure}[c]{0.18\textwidth}
        \includegraphics[width=1.0\textwidth, trim=600 100 500 100, clip]{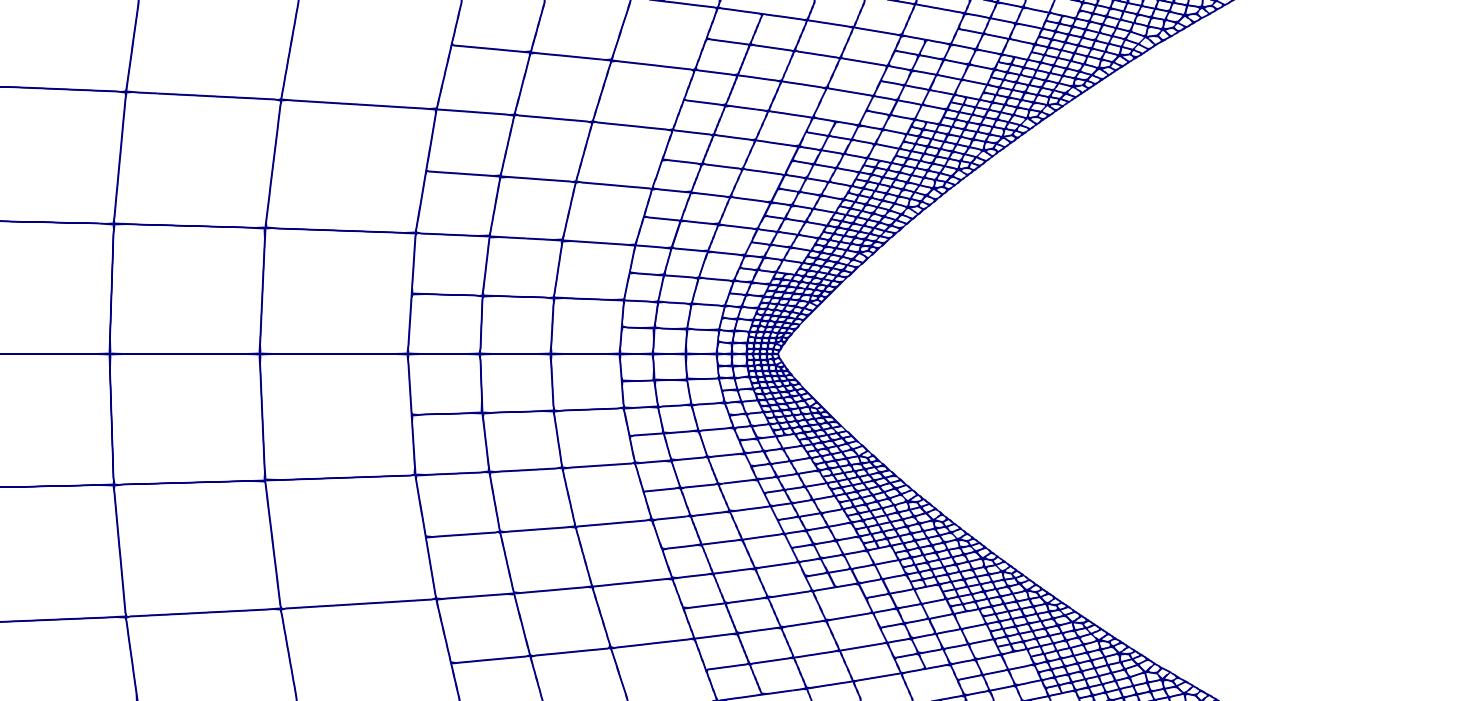}
        \subcaption{$p = 6$}
        \label{fig:cylinder2DMeshTipDetailP06}
    \end{subfigure}
    \caption{Influence of the $p$-value on the final grids in the vicinity of the upstream tip.}
    \label{fig:cylinder2DMeshTipDetails}
\end{figure*}

\subsubsection{Grid Deformation}
Besides the convergence, the attainable optimal shapes and the reduction of the objective function value, another major aspect refers to the quality of the mesh updates. Maintaining the grid quality during the optimization process is crucial to the success of the CAD-free optimization procedure. In this study, we focus on (a) the grid orthogonality near the boundaries as well as (b) the cell aspect ratio in the vicinity of the walls. It is seen, that using the $p$-Laplacian problem \eqref{eq:pLaplacianVariationalForm} with $p > 2$ significantly improves the quality of the mesh updates in comparison with updates obtained from  $p=2$. Figs.~\ref{fig:cylinder2DMeshTipDetailP02}~-~\ref{fig:cylinder2DMeshTipDetailP06} display the final grids in the upstream tip region for the five investigated values of $p$. A reasonable grid quality is generally observed for large values of $p$, even after substantial cumulative deformations due to several hundred optimization steps. On the contrary, the aspect ratio of the near wall cells in Fig.~\ref{fig:cylinder2DMeshTipDetailP02} increases significantly for $p=2$. The cells become stretched and tend to buckle in normal direction, which hampers the iterative convergence of the primal and adjoint flow solver. Therefore the procedure terminated after $354$ optimization steps for $p = 2$. In line with the change of the shape characteristics, a huge change of the mesh characteristics is observed when $p$ is increased from $p=2$ to $p=3$, cf.  Figs.~\ref{fig:cylinder2DMeshTipDetailP02} and ~\ref{fig:cylinder2DMeshTipDetailP03}. The grid is less compressed in the vicinity of the  upstream tip for $p=3$, despite the larger deviations from the initial grid indicated in Fig.~\ref{fig:compareFinalShapeP02P03Contours}. As outlined by Figs.~\ref{fig:cylinder2DMeshTipDetailP04}~-~\ref{fig:cylinder2DMeshTipDetailP06}, the grid does further improve when $p$ is augmented to $p = 4, 5, 6$. 
\begin{figure}[!ht]
    \centering
    \includegraphics[width=0.7\textwidth]{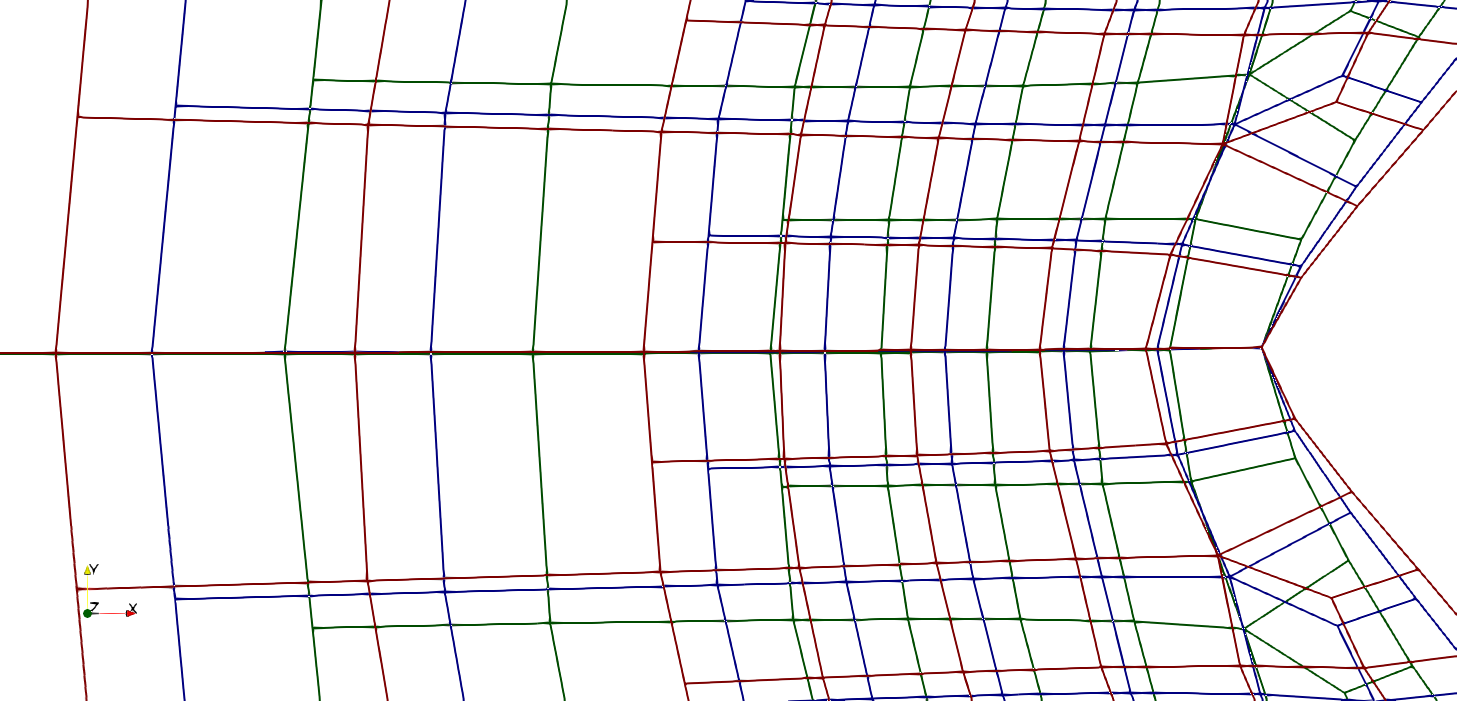}
    \caption{Superposition of upstream tip grids for $p = 4$ (green), $p = 5$ (blue) and $p = 6$ (red).}
    \label{fig:alignedGridsP04P05P06}
\end{figure}
A detailed comparison of the grids at the upstream tip of the final shape 
follows from Fig.~\ref{fig:alignedGridsP04P05P06} for $p = 4$ (green), $p = 5$ (blue) and for $p = 6$ (red). The post processed grids are aligned at the tip to support the comparison. The stretching of the grid increases with increasing $p$, which becomes obvious when observing the spacing along the horizontal center line.
\begin{figure*}[t]
    \centering
    \input{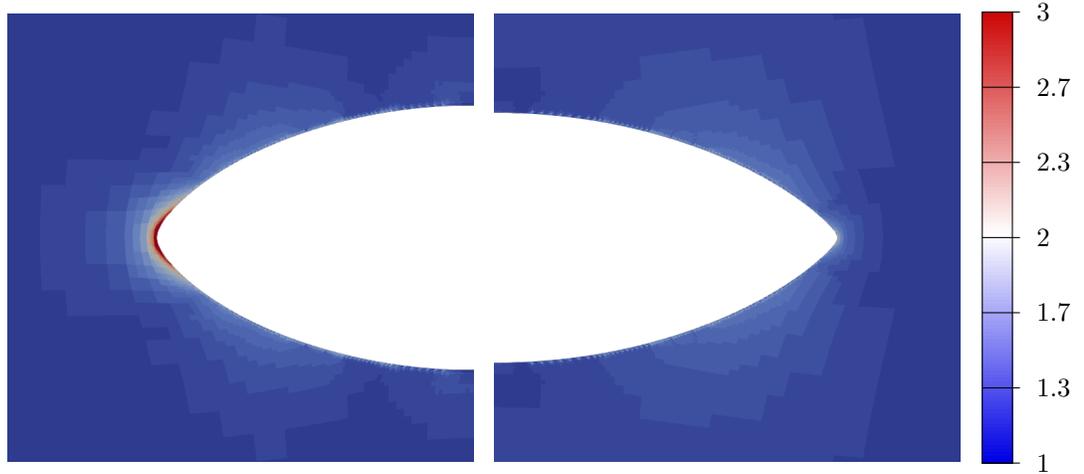}
    \caption{Cell aspect ratio of the final shapes obtained with $p = 2$ (left) and $p = 3$ (right).}
    \label{fig:cylinder2DMeshQualityAspectRatioP02P03}
\end{figure*}
Figure~\ref{fig:cylinder2DMeshQualityAspectRatioP02P03} depicts contour plots of the cell aspect ratio for the final grids obtained with $p = 2$ (left) and $p = 3$ (right). These $p$-values are particularly illustrative, since they denote a threshold for the characteristics of the mesh and the shape. The aspect ratio (AR) of the initial grid is generally close to $AR=1$. Only cells within the first two layers next to the cylinder boundary initially reach aspect ratios of up to $AR \le 1.5$. If one focuses on the final meshes, computations using $p = 2$ lead to a substantial amount of cells where the aspect ratio exceeds values of $AR=2.5$ and beyond, particularly  in the vicinity of the tip. The  maximum values approximately read $AR_{max}\le 11$. A significant improvement is achieved in conjunction with $p = 3$. Here, only few cells of the final mesh, located within a small area around the upstream tip, display aspect ratios greater than $AR>1.5$ and the peak values of the aspect ratio are limited to $AR_{max} \le 3$. 
\begin{figure*}[t]
    \centering
    \input{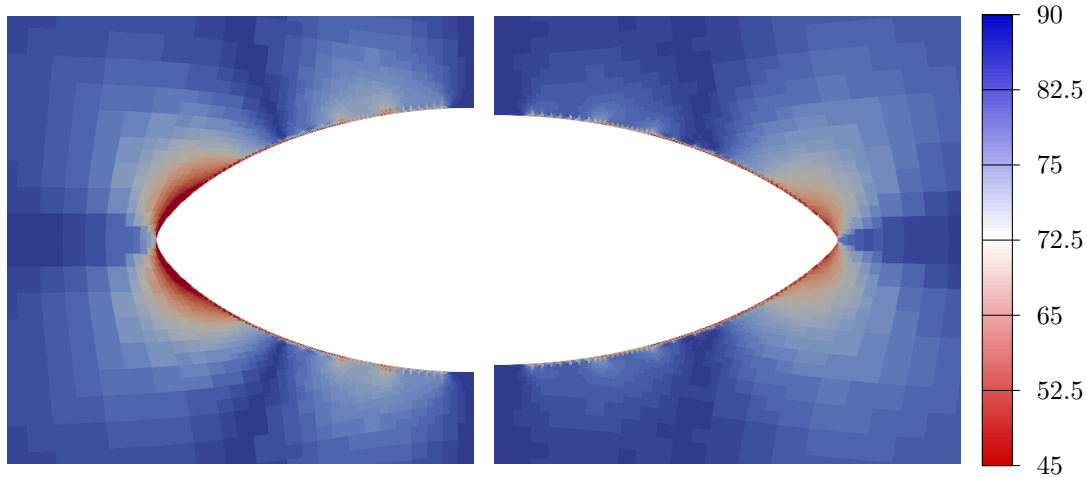}
    \caption{Minimum interior cell angle for the final meshes obtained with $p = 2$ (left) and $p = 3$ (right).}
    \label{fig:cylinder2DMeshQualityMinAnglep02P03}
\end{figure*}
The orthogonality of the cells within a mesh is of particular interest for the approximation of the boundary layer flow around the obstacle. Figure~\ref{fig:cylinder2DMeshQualityMinAnglep02P03} displays contour plots for the minimum interior angle of the final meshes using $p = 2$ (left) and $p = 3$ (right). Displayed angles reach from a less favorable value of $45^\circ$ (red) to a preferred value of $90^\circ$ (blue). Skewed cells occur close to tip and in regions where the shape tends towards a straight edge. This effect is much more pronounced for $p = 2$ (left) and significantly less obvious for  $p = 3$ (right). As illustrated by Fig. \ref{fig:alignedGridsP04P05P06}, the characteristic features of the mesh do not substantially change for higher $p$-values. 

\subsection{Drag Optimization in $3D$ Low Re Flow} 
\label{sc:application::3DLowRe}
\begin{figure*}[t]
    \centering
    \begin{subfigure}[c]{0.45\textwidth}
        \includegraphics[width=1.0\textwidth]{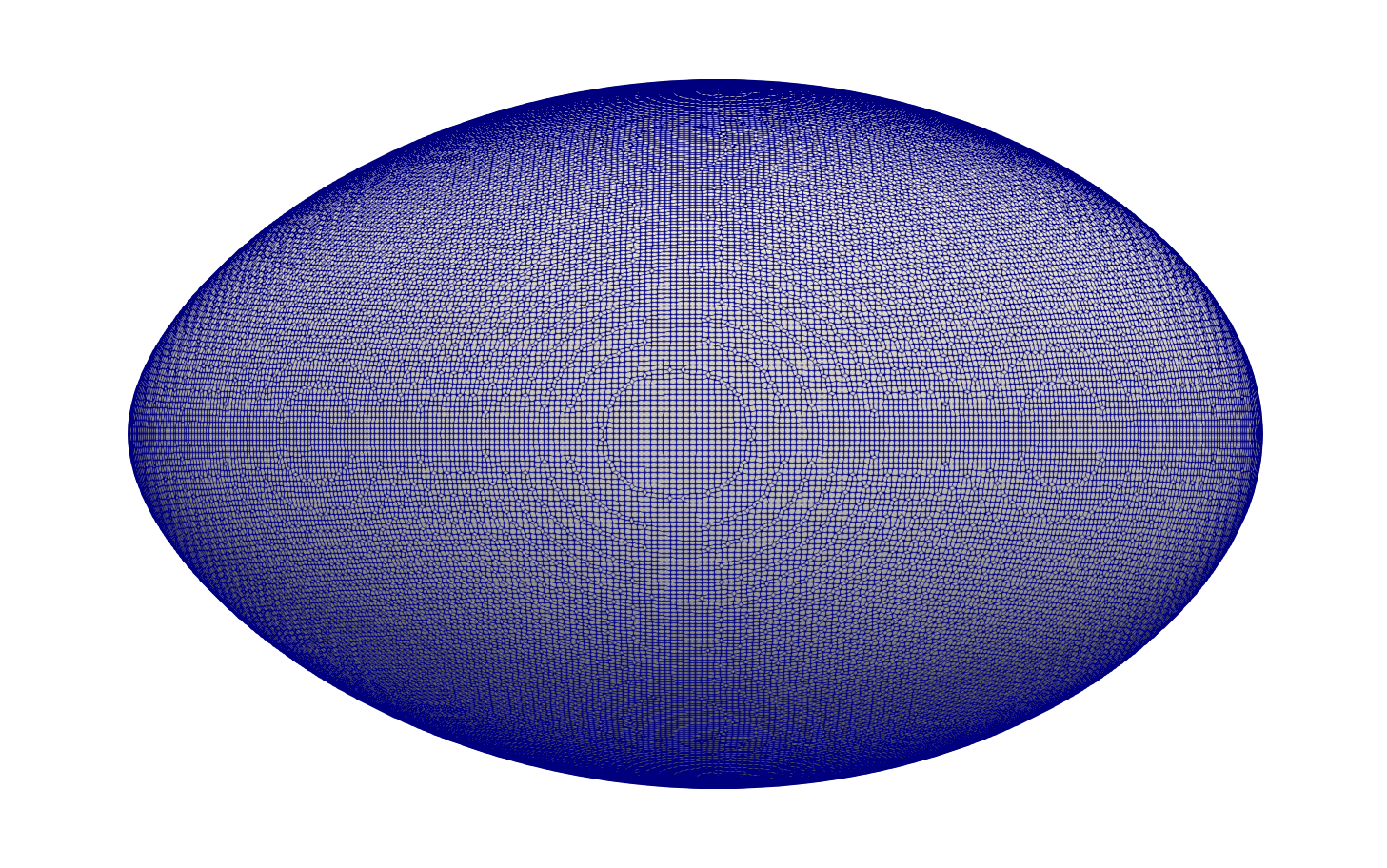}
        \subcaption{Side view.}
        \label{fig:sphere3DFinalShapep02:sideView}
    \end{subfigure}
    \begin{subfigure}[c]{0.45\textwidth}
        \includegraphics[width=1.0\textwidth]{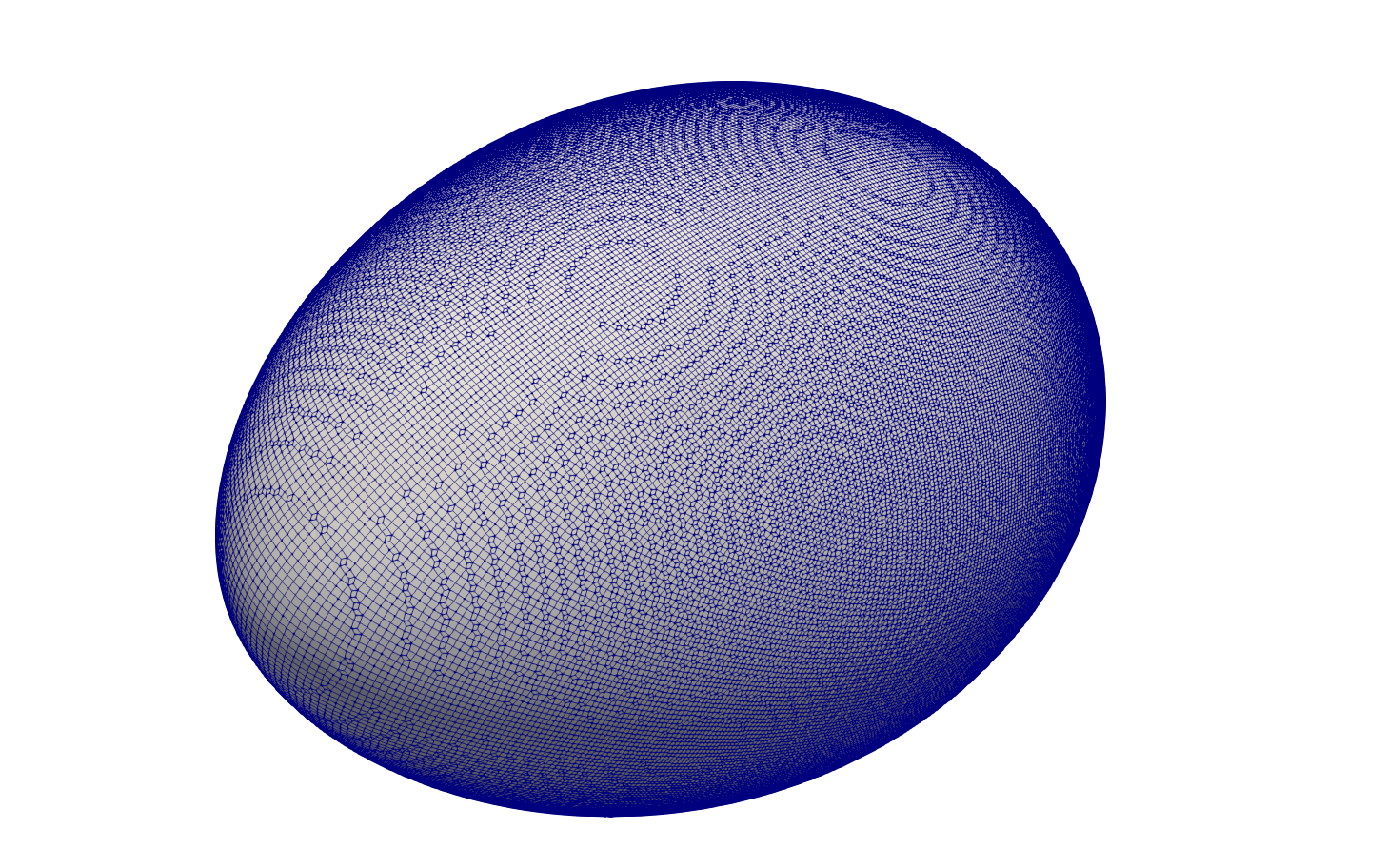}
        \subcaption{Perspective view onto the upstream tip.}
        \label{fig:sphere3DFinalShapep02:3dView}
    \end{subfigure}
    \caption{Drag optimization of a sphere at $Re=1$; final shape obtained from $p=2$ after $94$ steps.}
    \label{fig:sphere3DFinalShapep02}
\end{figure*}

This subsection is devoted to the drag optimization of a $3D$ unit-diameter sphere exposed to low Rey\-nolds number flow. Reported results are limited to the baseline case $p = 2$ and a single augmented $p$-level of $p=4$ -- which did already display substantial benefits in the $2D$ study. The computational domain, the position of the obstacle's barycenter and the boundary conditions agree with the $2D$ case outlined in section \ref{sc:application::2DLowRe}. The $2D$ domain is supplemented in  lateral direction by  $10 [m]$, and slip wall conditions are imposed along the lateral boundaries. The Rey\-nolds number compiled with the diameter and the approaching flow again reads $\mathit{Re} = 1$.
\begin{table}[!ht]
    \centering
    \begin{tabular}{c|c c c c}
        $p$   & Tip angle $[^\circ]$ & $\frac{1}{2}(\frac{a}{b} + \frac{a}{c}) _{(opt)}$ & $J/J^0$ & conv. step \\
        \hline
        $2$ &  $176.4$ & $1.61$ & $0.9507$ & $/$   \\
        $4$ &  $130.2$ & $2.13$ & $0.9360$ & $164$ \\
    \end{tabular}
    \caption{Performance indicators obtained with different $p$-values for the 3D low-Re case.}
    \label{tab:shpere3DoptimisationStats}
\end{table}
The boundary $\Gamma$ of the sphere is resolved by approximately $113$k wall adjacent cells and deforms under control. The total grid features about $1050$k control volumes. The wetted volume is restricted to conserve $\bar{c} = 5000 - \pi / 6 [m^3]$ and the initial parameters of the  augmented Lagrange algorithm~\ref{alg::augmented_lagrange} are again denoted in Tab.~\ref{tab:augmentedLagrangeInput}. The performance observed with the two investigated $p$-values is summarized in Tab.~\ref{tab:shpere3DoptimisationStats}. For $p = 4$ the optimization did converge after $164$ design steps and yields $6.4 \%$ drag reduction. Similar to the 2D case, the lower quality of the volume grid update restricted the amount design steps for $p=2$. Figure~\ref{fig:sphere3DFinalShapep02} depicts the last shape obtained from computations with $p = 2$ which referred to $94$ design steps and $4.93 \%$ drag reduction. Fig.~\ref{fig:sphere3DFinalShapep02:sideView} reveals that this shape is clearly characterized by round tips at the upstream and downstream ends.
\begin{figure*}[t]
    \centering
    \begin{subfigure}[c]{0.45\textwidth}
        \includegraphics[width=1.0\textwidth]{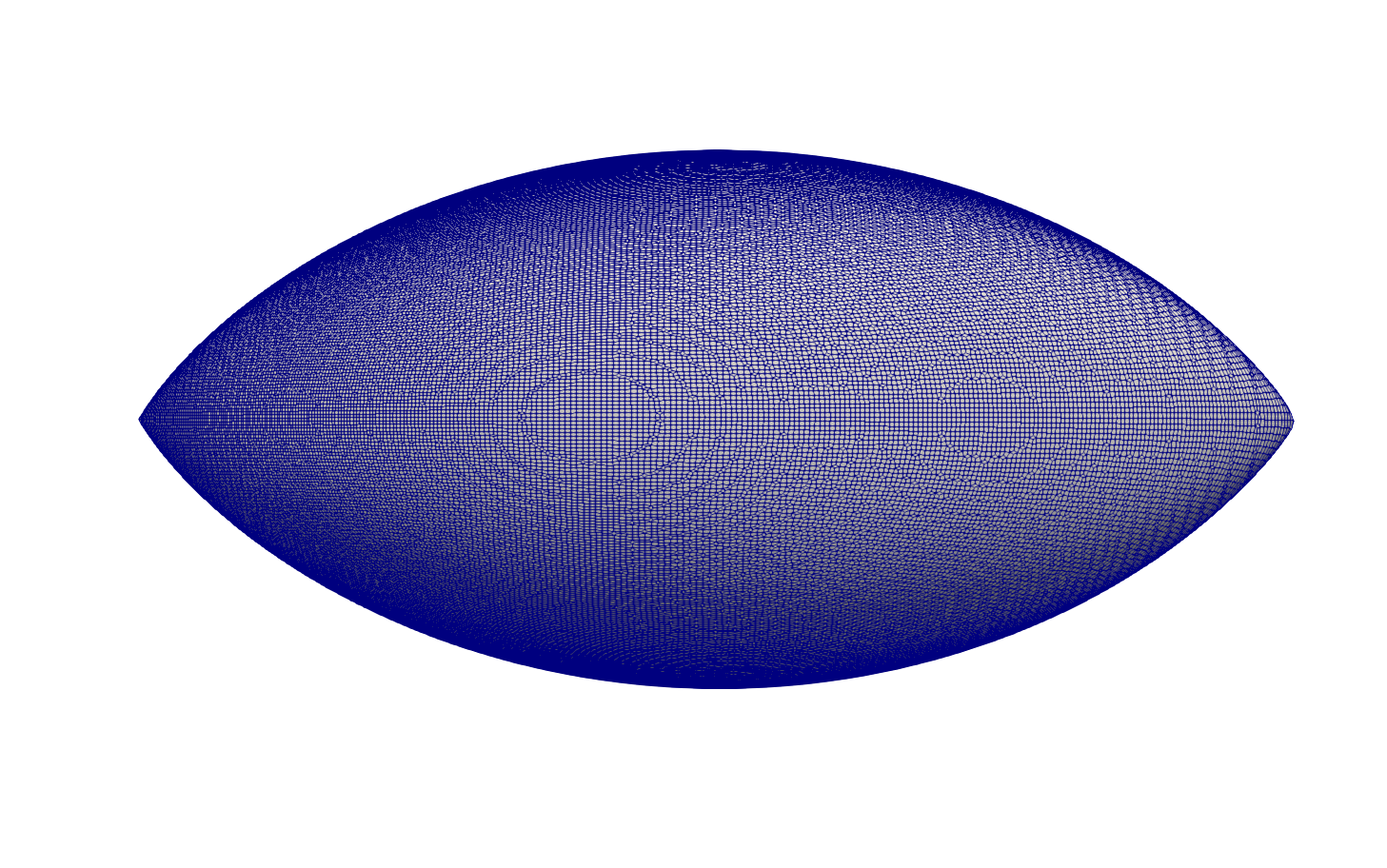}
        \subcaption{Side view.}
        \label{fig:sphere3DFinalShapep04:sideView}
    \end{subfigure}
    \begin{subfigure}[c]{0.45\textwidth}
        \includegraphics[width=1.0\textwidth]{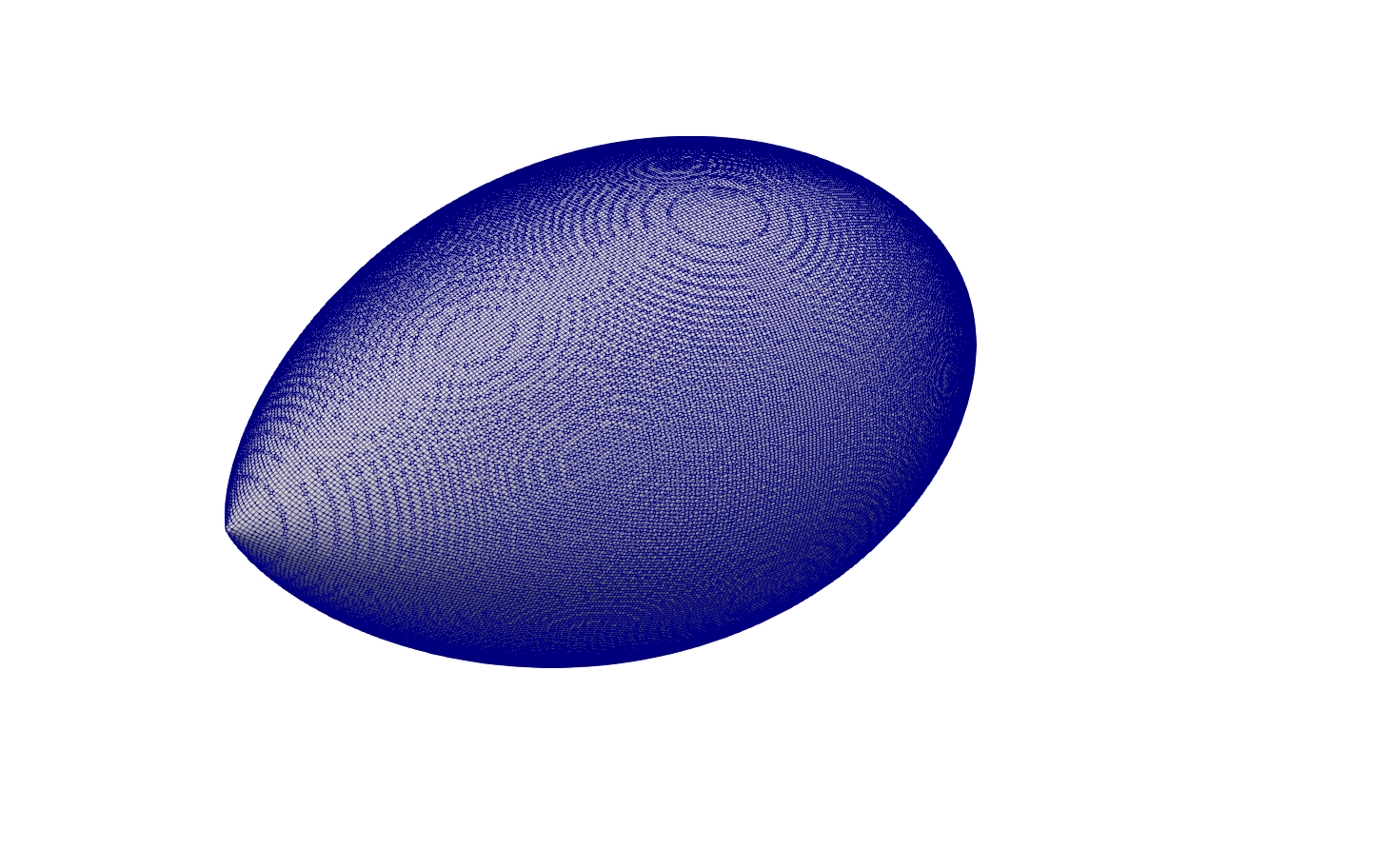}
        \subcaption{Perspective view onto the upstream tip.}
        \label{fig:sphere3DFinalShapep04:3dView}
    \end{subfigure}
    \caption{Drag optimization of a sphere at $Re=1$; optimal shape obtained from $p=4$ after $164$ steps.}
    \label{fig:sphere3DFinalShapeP04}
\end{figure*}
The optimal shape returned by $p = 4$ is shown in Fig.~\ref{fig:sphere3DFinalShapeP04}. As also indicated by the data listed in Tab.~\ref{tab:shpere3DoptimisationStats}, pointy upstream and downstream ends are seen for $p=4$ which also results in a significantly larger average half axis ratio.

\subsection{Drag Optimization in $2D$ Turbulent Flow}
\label{sc:application::2DHighReFlow}
Supplementary to the two low Reynolds number cases described above, we report the results obtained for a turbulent $2D$ drag optimization of an initial $a/b = 4/1$ ellipses where the longer half axis has a unit length of $a=1 [m]$. The incompressible fluid is characterized by a unit density and a dynamic viscosity of $\mu = 6.\overline{66} \cdot 10^{-7}$[Pa$\cdot$s]. The governing equations refer to RANS equations using a standard $k$-$\omega$ eddy-viscosity turbulence model \cite{wilcox1998turbulence} in combination with a wall-function approach. The Rey\-nolds number based on the longer axis of the ellipses and the approach flow reads $\mathit{Re} = 3 \cdot 10^6$. The computational domain and the boundary conditions agree with the information already used  in the first example, cf. section \ref{sc:application::2DLowRe}.  Results are again restricted to the baseline case $p=2$ and $p=4$. The wetted volume and the target barycenter are identical to those applied to the low Reynolds number experiment in section \ref{sc:application::2DLowRe}. The initial values for the algorithim~\ref{alg::augmented_lagrange} are denoted within the third colum of Tab.~\ref{tab:augmentedLagrangeInput}.
\begin{figure}[!ht]
    \centering
    \includegraphics[width=0.7\textwidth]{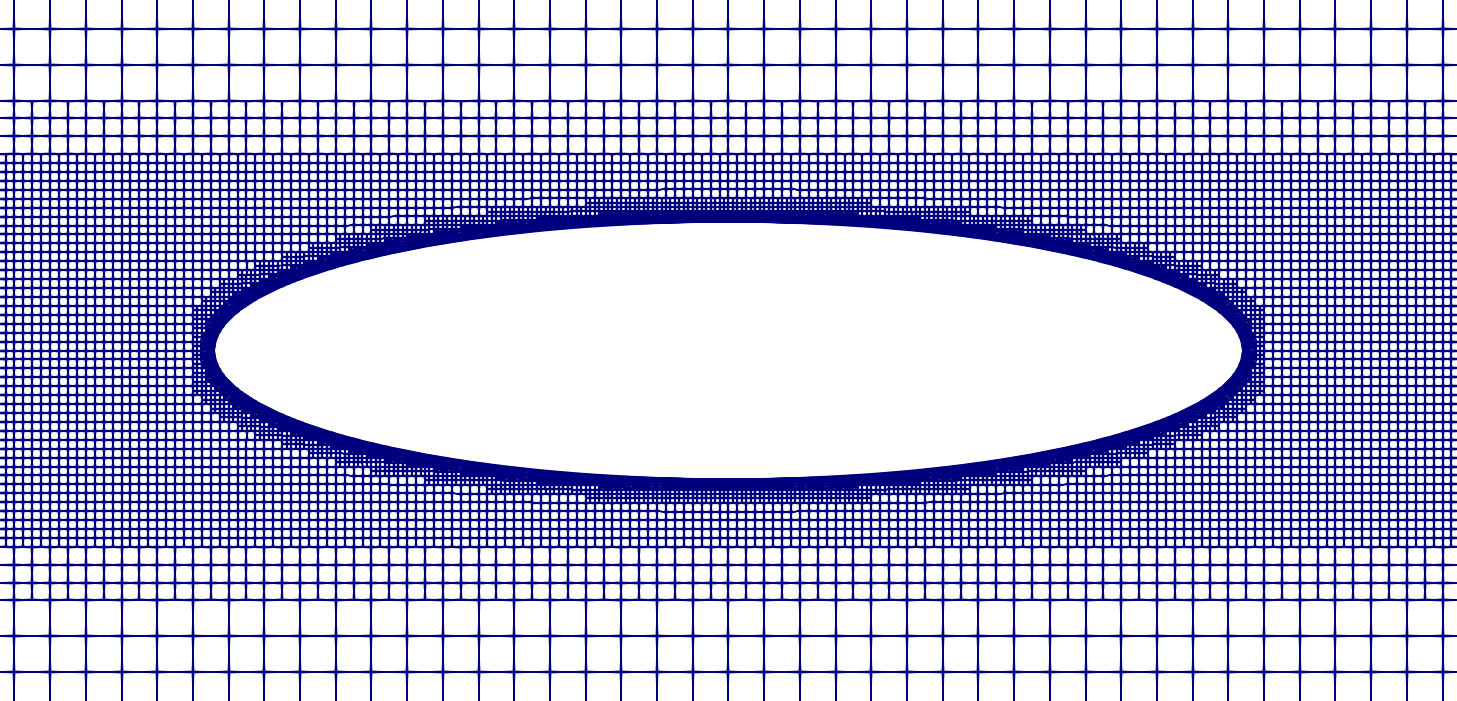}
    \caption{Initial mesh for the drag optimization of a $2D$ ellipses  exposed to horizontal approach flow at $Re=3 \cdot 10^6$.}
    \label{fig:ellipse2DInitialMesh}
\end{figure}
The initial grid again features local grid refinement near the boundary of the design surface to ensure an adequate resolution of the high curvature region and is depicted in Fig.~\ref{fig:ellipse2DInitialMesh}. The design surface $\Gamma$ is discretised by approximately $2300$ cells of equal size  and the volume grid features $46$k cells.
\begin{figure}[!ht]
    \centering
    \includegraphics[width=0.7\textwidth]{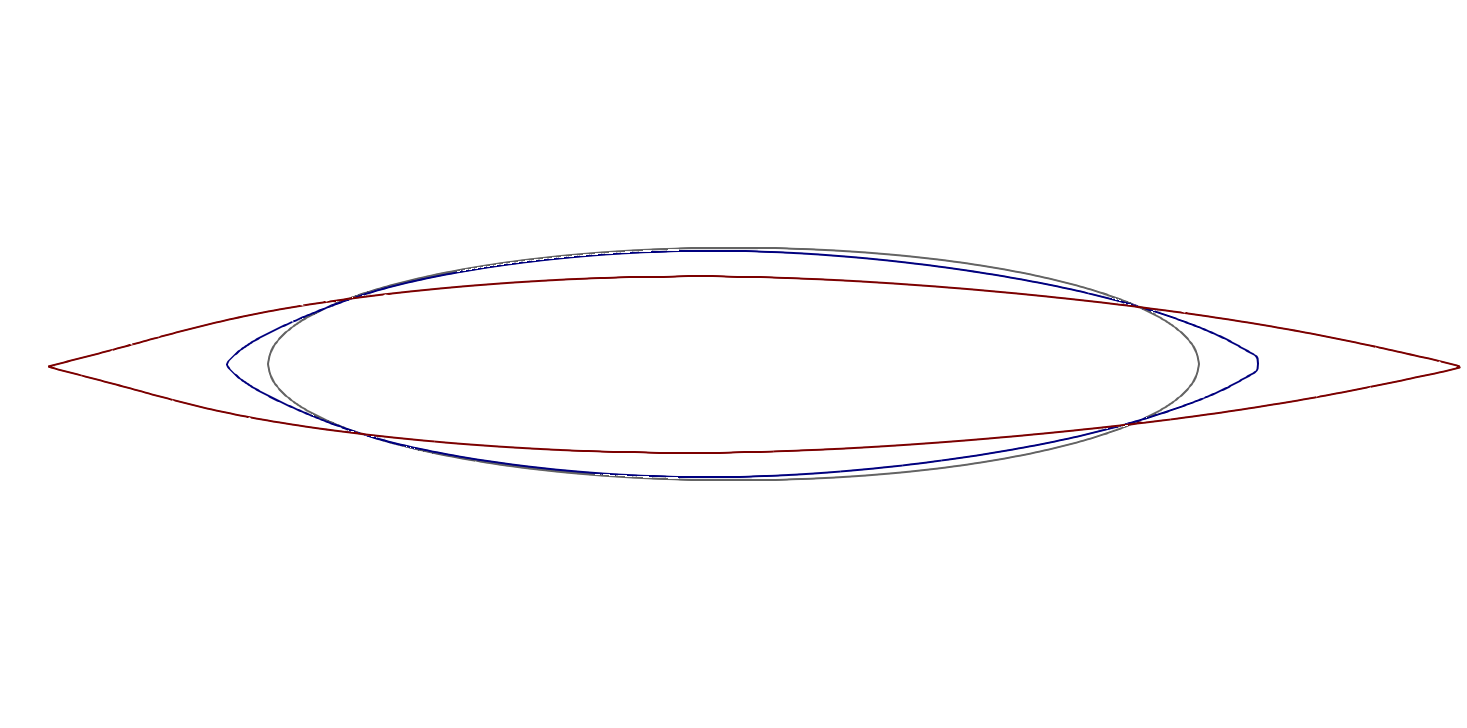}
    \caption{Drag optimization of a $2D$ elliptical cylinder exposed to horizontal approach flow at $Re=3 \cdot 10^6$; Comparison of the final shapes obtained for $p = 2$ (blue) and $p = 4$ (red) against the initial shape (grey).}
    \label{fig:ellipse2DContourCompare}
\end{figure}
The final shape contours obtained for $p = 2$ (blue) and  $p = 4$ (red) are outlined in Fig.~\ref{fig:ellipse2DContourCompare}. Data listed in Tab.~\ref{tab:ellipse2DoptimisationStats} reveals a significant difference of the shape and the drag reduction experienced with the different $p$-values. Similar to the $3D$ study, pointy upstream and downstream ends are seen for $p=4$. The larger $p$-value results in a significant increase of the initial half axis ratio while $a/b$ hardly increases for the smaller $p$-value. As in the previous studies, the shape optimization terminated before an optimal shape could be reached due to the severe distortion of the grid for the baseline value $p = 2$. For the baseline approach $p = 2$, the drag force of the body is reduced by $10.57 \%$  within $540$ optimization steps. However, for $p = 4$ the drag force is reduced by $32.6 \%$ when reaching the convergence criterion.
\begin{table}[!ht]
    \centering
    \begin{tabular}{c|c c c c}
        $p$   & Tip angle $[^\circ]$ & $a/b _{(opt)}$ & $J/J^0$ & conv. step \\
        \hline
        $2$ &  $149.5$ & $4.55$ & $0.894$ & $/$   \\
        $4$ &  $29.4$ & $8.02$ & $0.674$ & $589$ \\
    \end{tabular}
    \caption{Performance indicators obtained with different $p$-values for the $2D$ high Re case.}
    \label{tab:ellipse2DoptimisationStats}
\end{table}
\begin{figure*}[htp]
    \centering
    \begin{subfigure}[c]{0.4\textwidth}
        \includegraphics[width=1.0\textwidth]{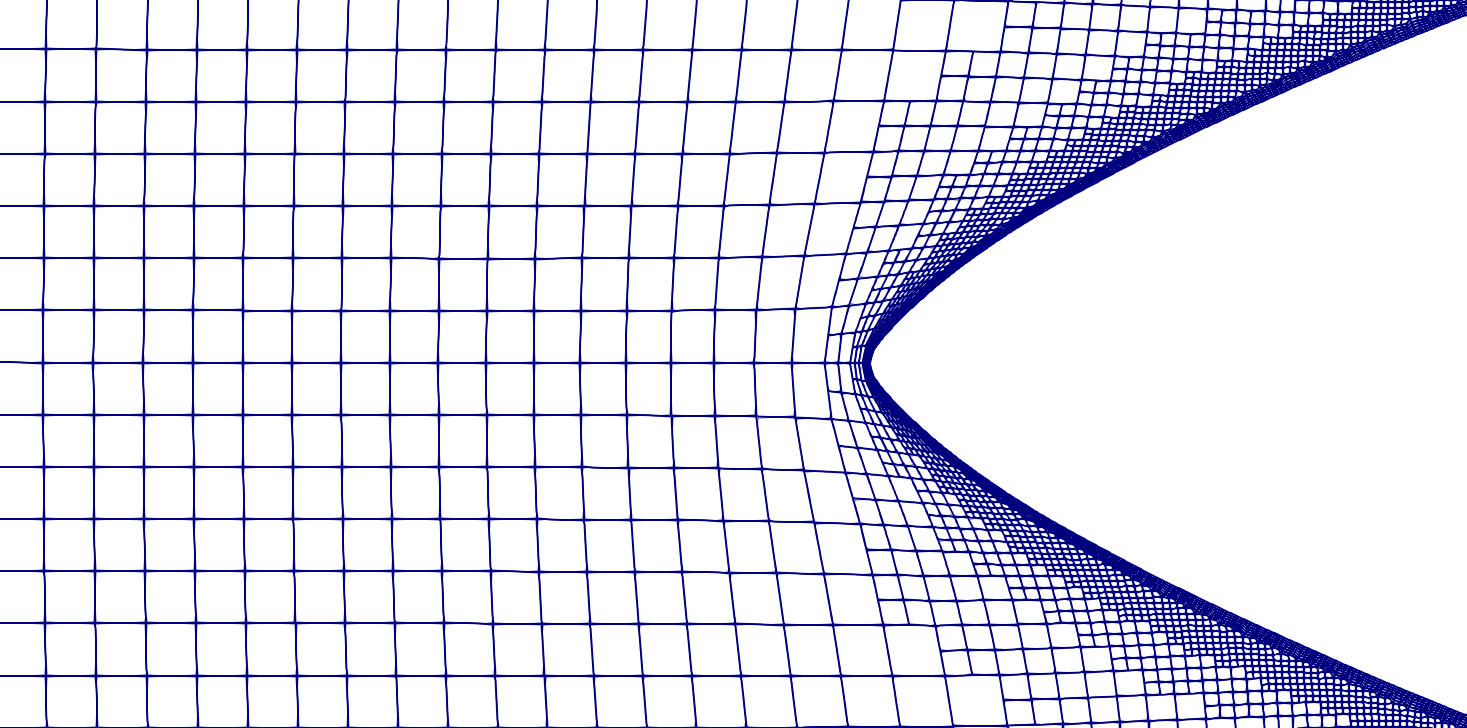}
        \subcaption{$p = 2$}
        \label{fig:ellipse2DMeshFrontTip:p02}
    \end{subfigure}
    \begin{subfigure}[c]{0.4\textwidth}
        \includegraphics[width=1.0\textwidth]{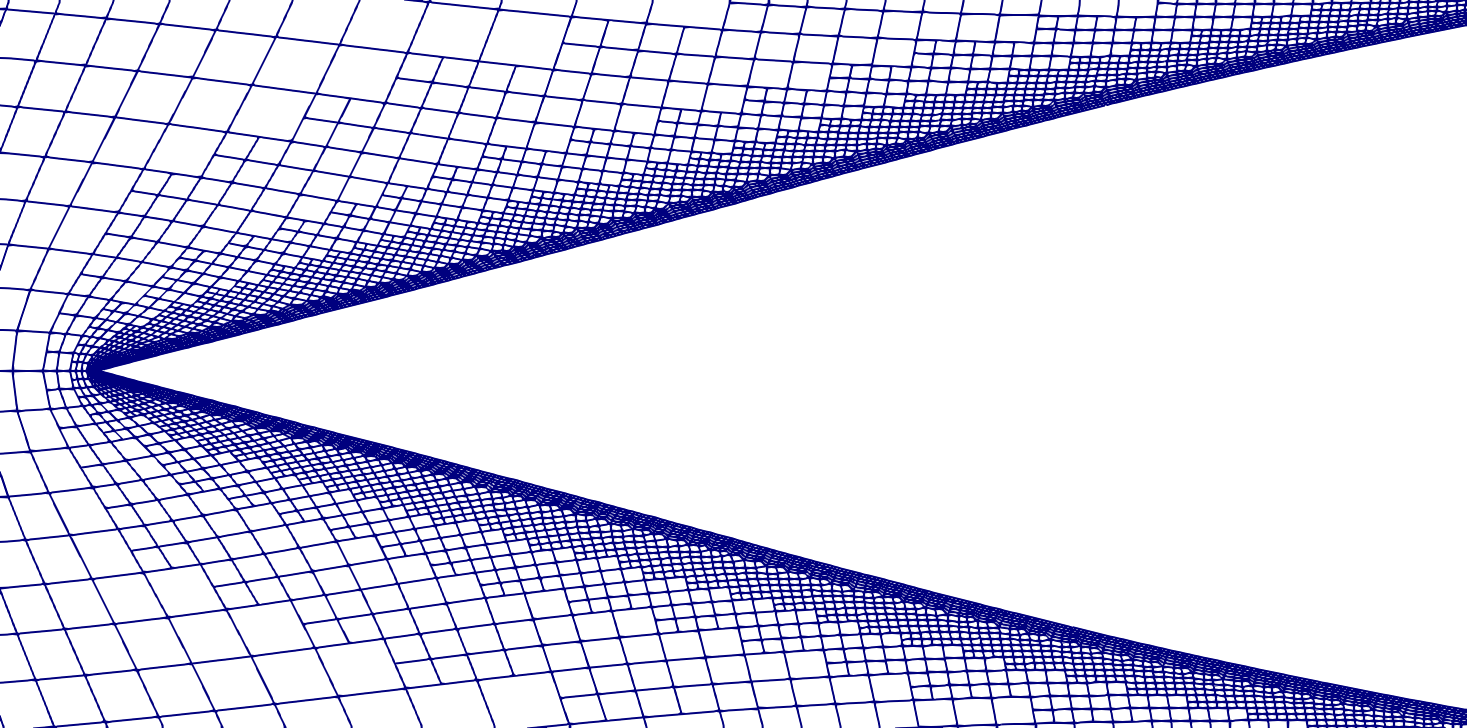}
        \subcaption{$p = 4$}
        \label{fig:ellipse2DMeshFrontTip:p04}
    \end{subfigure}
    \caption{Upstream tip of final shapes for the $2D$ high Reynolds case.}
    \label{fig:ellipse2DMeshFrontTip}
\end{figure*}
Figure~\ref{fig:ellipse2DMeshFrontTip} compares the final grids in the vicinity of the upstream tip for $p = 2$ (left) and $p = 4$ (right). As indicated by Fig.~\ref{fig:ellipse2DMeshFrontTip:p02}, the aspect ratio deteriorates for $p = 2$ since the near wall cells stretch in tangential direction. Moreover, cells (again) cluster at the tip. In contrast, the grid for $p = 4$ depicted in Fig.~\ref{fig:ellipse2DMeshFrontTip:p04} features evenly distributed cells along the design surface. The predicted shape displays a pointy upstream tip which is a much better approximation to the solution to the optimization problem.
\begin{figure*}[htp]
    \centering
    \begin{subfigure}[c]{0.4\textwidth}
        \includegraphics[width=1.0\textwidth]{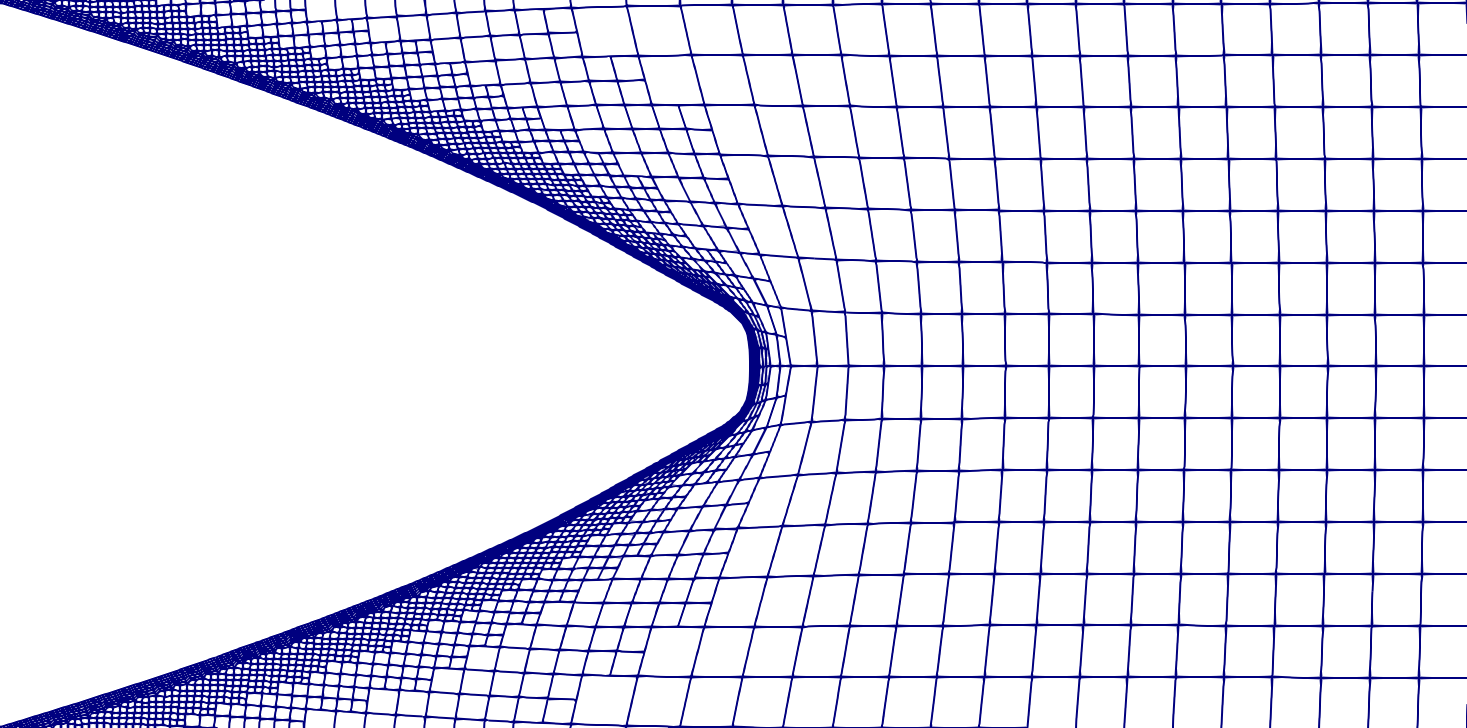}
        \subcaption{$p = 2$}
        \label{fig:ellipse2DMeshAftTip:p02}
    \end{subfigure}
    \begin{subfigure}[c]{0.4\textwidth}
        \includegraphics[width=1.0\textwidth]{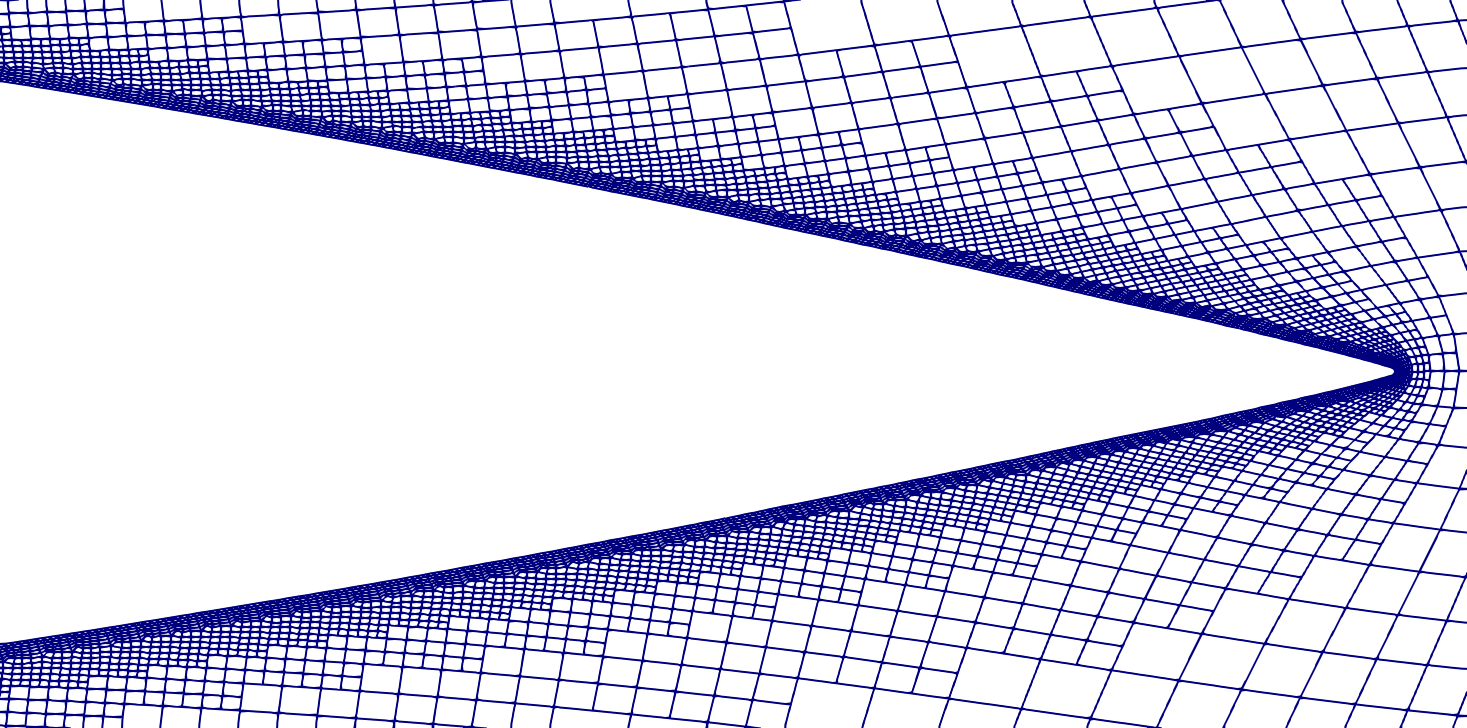}
        \subcaption{$p = 4$}
        \label{fig:ellipse2DMeshAftTip:p04}
    \end{subfigure}
    \caption{Downstream tip of final shapes for the 2D high Reynolds case.}
    \label{fig:ellipse2DMeshAftTip}
\end{figure*}
A similar conclusion follows from Fig.~\ref{fig:ellipse2DMeshAftTip} which describes the situation at the downstream end. It is observed, that the pointy ends develop at later stages of the optimization process, particularly in the rear. Thus, the (almost) pointy rear is hardly reached for $p = 2$. In combination with $p = 4$, the cell distribution at the rear remains evenly distributed. Mind that a small round tip can still observed in Fig.~\ref{fig:ellipse2DMeshAftTip:p04}. To receive a sharp pointy tip at the downstream end, a step size control would be beneficial.

\section{Conclusions}
\label{sc:conclusion}
We presented a novel approach for shape optimization by approximating Lipschitz continuous transformations based on the relaxation of the definition of the steepest descent direction. Examples included were restricted to fluid dynamic applications, but also  apply to other simulation areas. The main goal was to improve the shape optimization algorithm by considering a descent direction as the solution to the $p$-Laplace problem, and investigate the influence of increasing $p$ on the convergence of the shape optimization procedure as well as the obtained shapes and the related updates of the volume grid. An important aspect refers to the behavior of the relative differences  when $p$ is increased, since this might guide towards sufficiently high $p$ values associated with appreciated lower computational effort.

Results show that directions obtained from $p$-har\-monic solutions improve the convergence with increasing $p$. At the same time, the optimal shapes improve regarding the value of the objective function. A remarkable feature is related to the  ability of the $p$-Laplace approach  to yield shapes with edges or pointy shapes, even when the initial shape does not contain such features. Furthermore, the quality of the computational grid is virtually preserved even when large deformations of the initial shape occur and no  specific grid adjustment is required. Results of the present study suggest that $p=4$ seems a sufficiently large $p$-value to gather the benefits of the $p$-harmonic approach.

Within future research different solution algorithms for the $p$-Laplace problem may be considered to improve the approximation of the steepest descent direction and thus reduce the computational efforts. Moreover, applications to large-scale 3D problems may be investigated.

\section*{Acknowledments}
The authors acknowledge the support by the Deutsche Forschungsgemeinschaft (DFG) within the Research Training Group GRK 2583 "Modeling, Simulation and Optimization of Fluid Dynamic Applications” as well as within the research project ”Drag Optimisation of Ship Shapes” (Grant No. RU 1575/3-1). Michael Hinze acknowledges support of the DFG Priority Programme 1962 with projekt P8 "A Non-Smooth Phase-Field Approach to Shape Optimization with Instationary Fluid Flow". Selected computations were performed with resources provided by the North-German Super-computing Alliance (HLRN).

\section*{Authors' contributions}

\textbf{Peter Marvin M{\"u}ller}: Conceptualization, Methodology, Software, Validation, Formal analysis, Investigation, Writing - original draft, Writing - review \& editing., Visualization.
\textbf{Niklas K{\"u}hl}: Software, Validation, Formal analysis, Investigation, Writing - original draft, Writing - review \& editing
\textbf{Martin Siebenborn}: Funding acquisition, Conceptualization, Methodology, Writing - review \& editing.
\textbf{Klaus Deckelnick}: Idea, mathematical consulting.
\textbf{Michael Hinze}: Idea, mathematical consulting, methodology, Writing - review \& editing.
\textbf{Thomas Rung}: Project administration, Funding acquisition, Supervision, Conceptualization, Methodology, Resources, Writing - original draft, Writing - review \& editing.

\printbibliography

\end{document}